\documentclass[12pt,reqno]{amsart}

\usepackage[arrow,matrix,curve]{xy}

\usepackage[dvips]{graphicx}

\usepackage{amssymb, latexsym, amsmath, amscd, array, hyperref
}

\usepackage{setspace}

\theoremstyle{definition}

\numberwithin{equation}{section}

\newcommand\R {{\mathbb R}}
\newcommand\Q {{\mathbb Q}}

\newcommand\astr{{}^\ast\R}

\newcommand\megethos
{{$\mu\acute\varepsilon\gamma\varepsilon\vartheta{}o\varsigma$}}

\author[J.B.]{Jacques Bair}\address{J. Bair, HEC-ULG, University of
Liege, 4000 Belgium}\email{j.bair@ulg.ac.be}

\author[P.B.]{Piotr B\l{}aszczyk}\address{P. B\l{}aszczyk, Institute
of Mathematics, Pedagogical University of Cracow,
Poland}\email{pb@up.krakow.pl}

\author[R.E.]{Robert Ely}\address{R. Ely, Department of Mathematics,
University of Idaho, Moscow, ID 83844 US}\email{ely@uidaho.edu}

\author[V.H.]{Val\'erie Henry}\address{V.~Henry, Department of
Mathematics, University of Namur, 5000
Belgium}\email{vhen@math.fundp.ac.be}

\author[V.K.]{Vladimir Kanovei} \address{V. Kanovei, IPPI, Moscow, and
MIIT, Moscow, Russia} \email{kanovei@rambler.ru}

\author[K.K.]{Karin U. Katz}\address{K. Katz, Department of
Mathematics, Bar Ilan University, Ramat Gan 52900 Israel}

\author[M.K.]{Mikhail G. Katz}\address{M. Katz, Department of
Mathematics, Bar Ilan University, Ramat Gan 52900 Israel}
\email{katzmik@macs.biu.ac.il}

\author[S.K.]{Semen S. Kutateladze}\address{S. Kutateladze, Sobolev
Institute of Mathematics, Novosibirsk State University, Russia}
\email{sskut@math.nsc.ru}

\author[T.M.]{Thomas McGaffey}\address{T. McGaffey, Rice University,
US}\email{thomasmcgaffey@sbcglobal.net}

\author[P.R]{Patrick Reeder}\address{P. Reeder, Holyoke College}
\email{patreeder@gmail.com}

\author[D.S.]{David M. Schaps}\address{D. Schaps, Department of
Classical Studies, Bar Ilan University, Ramat Gan 52900 Israel}
\email{dschaps@mail.biu.ac.il}

\author[D.S.]{David Sherry}\address{D. Sherry, Department of
Philosophy, Northern Arizona University, Flagstaff, AZ 86011,
US}\email{david.sherry@nau.edu}

\author[S.S.]{Steven Shnider}\address{S.~Shnider, Department of
Mathematics, Bar Ilan University, Ramat Gan 52900
Israel}\email{shnider@macs.biu.ac.il}

\begin{document}


\title[Interpreting Leibniz and Euler]{Interpreting the infinitesimal
mathematics of Leibniz and Euler}

\subjclass[2000]{Primary 01A50; 
Secondary
26E35,       
01A85,       
03A05%
}

\begin{abstract}
We apply Benacerraf's distinction between mathematical \emph{ontology}
and mathematical \emph{practice} (or the \emph{structures}
mathematicians use in practice) to examine contrasting interpretations
of infinitesimal mathematics of the 17th and 18th century, in the work
of Bos, Ferraro, Laugwitz, and others.  We detect Weierstrass's ghost
behind some of the received historiography on Euler's infinitesimal
mathematics, as when Ferraro proposes to \emph{understand} Euler in
terms of a Weierstrassian notion of limit and Fraser declares
classical analysis to be a ``primary point of reference for
understanding the eighteenth-century theories''.  Meanwhile, scholars
like Bos and Laugwitz seek to explore Eulerian \emph{methodology},
\emph{practice}, and \emph{procedures} in a way more faithful to
Euler's own.

Euler's use of infinite integers and the associated infinite products
is analyzed in the context of his infinite product decomposition for
the sine function.  Euler's \emph{principle of cancellation} is
compared to the Leibnizian \emph{transcendental law of homogeneity}.
The Leibnizian \emph{law of continuity} similarly finds echoes in
Euler.  

We argue that Ferraro's assumption that Euler worked with a
\emph{classical} notion of quantity is symptomatic of a
post-Weierstrassian placement of Euler in the Archimedean track for
the development of analysis, as well as a blurring of the distinction
between the dual tracks noted by Bos.  Interpreting Euler in an
Archimedean conceptual framework obscures important aspects of Euler's
work.  Such a framework is profitably replaced by a syntactically more
versatile modern infinitesimal framework that provides better proxies
for his inferential moves.

Keywords: Archimedean axiom; infinite product; infinitesimal; law of
continuity; law of homogeneity; principle of cancellation; procedure;
standard part principle; ontology; mathematical practice; Euler;
Leibniz

\end{abstract}

\maketitle

\noindent

\tableofcontents

\section{Introduction}
\label{s1}

This text is part of a broader project of re-appraisal of the
Leibniz--Euler--Cauchy tradition in infinitesimal mathematics that
Weierstrass and his followers broke with around 1870.

In the case of Cauchy, our task is made easier by largely
traditionalist scholars G.~Schubring and G.~Ferraro.  Thus, Schubring
distanced himself from the Boyer--Grabiner line on Cauchy as the one
\emph{who gave you the epsilon} in the following terms: ``I am
criticizing historiographical approaches like that of Judith Grabiner
where one sees epsilon-delta already realized in Cauchy''
\cite[Section~3]{Sc15}.  Ferraro goes even further and declares:
``Cauchy uses infinitesimal neighborhoods of~$x$ in a decisive
way\ldots{} Infinitesimals are not thought as a mere \emph{fa\c con de
parler}, but they are conceived as numbers, though a theory of
infinitesimal numbers is lacking.''  \cite[p.~354]{Fe08} Ferraro's
comment is remarkable for two reasons:
\begin{itemize}
\item
it displays a clear grasp of the procedure \emph{vs} ontology
distinction (see below Section~\ref{11b});
\item
it is a striking recognition of the \emph{bona fide} nature of
Cauchy's infinitesimals that is a clear break with Boyer--Grabiner.
\end{itemize}
Ferraro's comment is influenced by Laugwitz's perceptive analysis of
Cauchy's sum theorem in \cite{La87}, a paper cited several times on
Ferraro's page 354.  For further details on Cauchy see the articles
\cite{BKS}, \cite{Bl16}.

In this article, we propose a re-evaluation of Euler's and, to an
extent, Leibniz's work in analysis.  We will present our argument in
four stages of increasing degree of controversy, so that readers may
benefit from the text even if they don't agree with all of its
conclusions.

(1) We argue that Euler's procedures in analysis are best proxified in
modern infinitesimal frameworks rather than in the received modern
Archimedean ones, by showing how important aspects of his work have
been underappreciated or even denigrated because inappropriate
conceptual frameworks are being applied to interpret his work.  To
appreciate properly Euler's work, one needs to abandon extraneous
\emph{ontological} matters such as the continuum being punctiform
(i.e., made out of points) or nonpunctiform, and focus on the
\emph{procedural} issues of Euler's actual mathematical practice.

(2) One underappreciated aspect of Euler's work in analysis is its
affinity to Leibniz's.  A number of Eulerian procedures are consonant
with those found in Leibniz, such as the law of continuity (governing
the passage from an Archimedean continuum to an infinitesimal-enriched
continuum) and the transcendental law of homogeneity (governing the
passage from an infinitesimal-enriched continuum back to an
Archimedean continuum).  This is consistent with the teacher--student
lineage from Leibniz to Johann Bernoulli to Euler.

(3) Leibniz wrote in 1695 that his infinitesimals violate the property
expressed by Euclid's Definition V.4 (see \cite[p.\;288]{Le95a}).%
\footnote{Actually Leibniz referred to V.5; in some editions of the
\emph{Elements} this Definition does appear as V.5.  Thus, \cite{Eu60}
as translated by Barrow in 1660 provides the following definition in
V.V (the notation ``V.V'' is from Barrow's translation):
\emph{Those numbers are said to have a ratio betwixt them, which being
multiplied may exceed one the other}.}
This axiom is a variant of what is known today as the Archimedean
property.  Thus, Leibnizian infinitesimals violate the Archimedean
property when compared to other quantities.

(4) Our reading is at odds with the \emph{syncategorematic}
interpretation elaborated in \cite[Chapter~5]{Is}, \cite{Ar08}, and
elsewhere.  Ishiguro, Arthur, and others maintain that Leibniz's
continuum was Archimedean, and that his infinitesimals do not
\emph{designate} and are \emph{logical fictions} in the sense of
Russell.  The leap by Ishiguro (and her followers) from infinitesimals
being \emph{fictions} to their being \emph{logical fictions} is a
non-sequitur analyzed in articles in \emph{Erkenntnis} \cite{KS1} and
in \emph{Studia Leibnitiana} \cite{SK}.  Arthur's interpretation was
also challenged in \cite{Th12}.  The fictions in question are
\emph{pure} rather than \emph{logical}, meaning that they do
\emph{designate} insofar as our symbolism allows us to think about
infinitesimals.  This is consistent with interpretations of Leibniz by
\cite{Bos} and \cite{Je15} (see Section~\ref{sync1} for a discussion
of Jesseph's analysis).  Euler similarly works explicitly with
infinite and infinitesimal numbers rather than some kind of paraphrase
thereof in terms of proto-Weierstrassian hidden quantifiers.

In Appendix~\ref{31}, we examine the mathematical details of the
Eulerian procedures in the context of his proof of the infinite
product decomposition for the sine function and related results.

In addition to Robinson's framework, other modern theories of
infinitesimals are also available as possible frameworks for the
interpretation of Euler's procedures, such as Synthetic Differential
Geometry \cite{Ko06}, \cite{Be08} and Internal Set Theory \cite{Ne77},
\cite{kr}.  See also \cite{NK} as well as \cite{KKN}.  Previous
studies of the history of infinitesimal mathematics include
\cite{KK11}, \cite{BK}, \cite{Ba13}, \cite{KSS13}, \cite{Ca13},
\cite{Ba14}, \cite{Ka15a}.

\section{Historiography}
\label{one}

It is a subject of contention among scholars whether science
(including mathematics) develops continuously or by discontinuous
leaps.  The use of paradigm shifts by \cite{Ku} is the most famous
instance of the discontinuous approach.  The discontinuous case is
harder to make for mathematics than for the physical sciences: we gave
up on phlogiston and caloric theory, but we still use the Pythagorean
theorem and l'H\^opital's rule.

\subsection{Continuity and discontinuity}
\label{s21}

We argue that the continuous \emph{vs} discontinuous dichotomy is
relevant to understanding some of the current debates in interpreting
classical infinitesimalists like Leibniz and Euler.  Thus, A.~Robinson
argued for continuity between the Leibnizian framework and his own,
while H.~Bos rejected Robinson's contention in the following terms:
\begin{quote}
\ldots{} the most essential part of non-standard analysis, namely the
proof of the existence of the entities it deals with, was entirely
absent in the Leibnizian infinitesimal analysis, and this constitutes,
in my view, so fundamental a difference between the theories that the
Leibnizian analysis cannot be called an early form, or a precursor, of
non-standard analysis \cite[p.\;83]{Bos}.
\end{quote}

Of course, many scholars reject continuity not merely between
Robinson's framework and historical infinitesimals, but also between
the received modern mathematical frameworks and historical
infinitesimals.  A case in point is Ferraro's treatment of an
infinitesimal calculation found in \cite[pp.\;11-12]{Eu30}.  Here
Euler sought the value of the ratio
$\big({1-x^{g/(f+g)}_{\phantom{I}}}\big)/{g}$ for~$f = 1$ and~$g=0$ by
applying l'H\^opital's rule to~$\frac{1-x^z}{z}$.  Ferraro proceeds to
present the problem ``from a modern perspective'' by analyzing the
function~$f(z)=\frac{1-x^z}{z}$ and its behavior near~$z=0$ in the
following terms:

\begin{quote}
From the modern perspective, the problem of extending the function
$f(z)=\frac{1-x^z}{z}$ in a continuous way means that\ldots{} the
domain~$D$ of~$f(z)$ has a point of accumulation at 0 so that we can
attempt to calculate the limit as~$z\to 0$, where
by~$\lambda=\lim_{z\to c}f(z)$ [the~$c$ in Ferraro's formula needs to
be replaced by~$0$] we mean: given any~$\varepsilon>0$ there exists
a~$\delta>0$ such that if~$z$ belongs to~$D$ and~$|z|<\delta$
then~$|f(z)-\lambda|<\varepsilon$; \ldots{} This procedure is
\emph{substantially meaningless} for Euler
\cite[p.\;46]{Fe04}. (emphasis added)
\end{quote}
Ferraro's concluding remarks concerning ``substantially meaningless''
procedures place him in the discontinuity camp.

While there is a great deal of truth in the discontinuous position,
particularly with regard to currently prevalent ontological frameworks
(set-theoretic or category-theoretic), we will argue for a limited
reading of the history of analysis from the perspective of
\emph{continuous} development in the following sense.  As we analyze
the history of analysis since the 17th century, we note stark
differences among the objects with which mathematicians reason; there
are for example no \emph{sets} as explicit mathematical objects in
Leibniz or Euler.  On the other hand, there are important continuities
in the principles which guide the inferences that they draw; for
example, Leibniz's \emph{transcendental law of homogeneity}, Euler's
\emph{principle of cancellation}, and the \emph{standard part
principle} exploited in analysis over a hyperreal
extension~$\mathbb{R}\subseteq{}^\ast\mathbb{R}$.

The crucial distinction here is between practice and ontology, as we
detail below in Section~\ref{11b}.  We will argue that there is a
historical continuity in mathematical practice but discontinuity in
mathematical ontology.  More specifically, the set-theoretic semantics
that currently holds sway is a discontinuity with respect to the
historical evolution of mathematics.  Scholars at times acknowledge
the distinction in relation to their own work, as when Ferraro speaks
about the \emph{intensional} nature of the entities in Euler in
\cite[p.\;44]{Fe04} and the \emph{syntactic} nature of algebraic and
analytic operations \cite[p.\;203]{Fe08}, but not always when it comes
to passing judgment on Laugwitz's work; see Section~\ref{s22}.

\subsection{Procedures and proxies}
\label{s22}

In the case of Euler, we will examine philosophical issues of
interpretation of infinitesimal mathematics (more specifically, the
use of infinitesimals and infinite integers) and seek to explore the
roots of the current situation in Euler scholarship, which seems to be
something of a dialog of the deaf between competing approaches.  Some
aspects of Euler's work in analysis were formalized in terms of modern
infinitesimal theories by Laugwitz, McKinzie, Tuckey, and others.
Referring to the latter, G. Ferraro claims that ``one can see in
operation in their writings a conception of mathematics which is
\emph{quite extraneous} to that of Euler'' \cite[p.\;51]{Fe04}
(emphasis added).  Ferraro concludes that ``the attempt to specify
Euler's notions by applying modern concepts is only possible if
elements are used which are \emph{essentially alien} to them, and thus
Eulerian mathematics is transformed into something wholly different.''
\cite[pp.\;51-52]{Fe04} (emphasis added)

Now \emph{quite extraneous} and \emph{essentially alien} are strong
criticisms.  The vagueness of the phrase ``to specify Euler's notions
by applying modern concepts" makes it difficult to evaluate Ferraro's
claim here.  If \emph{specification} amounts to bringing to light
tacit assumptions in Euler's reasoning, then it is hard to see why
Ferraro uses such harsh language.

We find a different attitude in P.~Reeder's approach to Euler.  Reeder
writes:
\begin{quote}
I aim to reformulate a pair of proofs from [Euler's]
\emph{Introductio} using concepts and techniques from Abraham
Robinson's celebrated non-standard analysis (NSA).  I will
specifically examine Euler's proof of the Euler formula and his proof
of the divergence of the harmonic series.  Both of these results have
been proved in subsequent centuries using epsilontic (standard
epsilon-delta) arguments.  The epsilontic arguments differ
significantly from Euler's original proofs.  \cite[p.\;6]{Re12}.
\end{quote}
Reeder concludes that ``NSA possesses the tools to provide appropriate
proxies of the inferential moves found in the \emph{Introductio}."
Reeder finds significant similarities between some of Euler's proofs
and proofs in a hyperreal framework.  Such similarities are missing
when one compares Euler's proofs to proofs in the~$\epsilon,\delta$
tradition.  We take this to mean that Euler's conception has more in
common with the syntactic resources available in a modern
infinitesimal tradition than in the~$\epsilon,\delta$ tradition.

Scholars thus appear to disagree sharply as to the relevance of modern
theories to Euler's mathematics, and as to the possibility of
meaningfully reformulating Euler's infinitesimal mathematics in terms
of modern theories.

\subsection{Precalculus or analysis?}

Having mentioned Euler's \emph{Introductio}, we would like to clarify
a point concerning the nature of this book.  Blanton writes in his
introduction that ``the work is strictly pre-calculus''
\cite[p.\;xii]{aie}.  Is this an accurate description of the book?  It
is worth keeping the following points in mind.
\begin{enumerate}
\item The algebraic nature of the \emph{Introductio} was mirrored 70
years later by Cauchy's \emph{Cours d'Analyse}, which was subtitled
\emph{Analyse Algebrique}.  Laugwitz noted in fact that \emph{Cours
d'Analyse} was modeled on Euler's \emph{Introductio}
\cite[p.\;52]{La99}.
\item There may not be much material related to differentiation in
\emph{Introductio}, but series are dealt with extensively.  Series
certainly being part of analysis, it seems more reasonable to describe
\emph{Introductio} as \emph{analysis} than precalculus.
\item 
\label{i3}
Infinitesimals in \emph{Introductio} and differentials in
\emph{Institutiones} are arguably of similar nature.  Leibniz already
thought of differentials as infinitesimals, as did Johann Bernoulli.
There is little reason to assume otherwise as far as Euler is
concerned, particularly since he viewed all his analysis books as a
unified whole.
\end{enumerate}
To elaborate further on item~\eqref{i3}, note that Euler writes in his
\emph{Institutionum calculi integralis} as follows: ``In calculo
differentiali iam notavi, quaestionem de differentialibus non absolute
sed relative esse intelligendam, ita ut, si~$y$ fuerit function
quaecunque ipsius~$x$, non tam ipsum eius differentiale~$dy$, quam
eius ratio ad differentiale~$dx$ sit definienda.''  \cite[p.\;6,
Scholion~1, \S\,5]{Eu68} This can be translated as follows: ``Now in
differential calculus I have observed that an investigation of
differentiation is to be understood as not absolute but relative;
namely, if~$y$ is a function of~$x$, what one needs to define is not
so much its differential~$dy$ itself as its ratio to the
differential~$dx$.''

The comment indicates that throughout the period 1748-1768, Euler
thinks of infinitesimals and differentials as essentially
interchangeable.

\subsection
{Practice versus ontology}
\label{11b}

In an influential essay ``The Relation Between Philosophy of Science
and History of Science,'' M. Wartofsky argues that historiography of
science needs to begin its analysis by mapping out an ontology of the
scientific field under investigation.  Here ontology is to be
understood in a broader sense than merely the ontology of the entities
exploited in that particular science -- such as numbers, functions,
sets, etc., in the case of mathematcs -- but rather develop the
ontology of mathematics as a scientific theory itself
\cite[p.\;723]{Wa76}.  

As a modest step in this direction we distinguish between the
(historically relative) ontology of the mathematical objects in a
certain historical setting, and its procedures, particularly
emphasizing the different roles these components play in the history
of mathematics.  More precisely, our \emph{procedures} serve as a
representative of what Wartofsky called the \emph{praxis}
characteristic of the mathematics of a certain time period, and our
\emph{ontology} takes care of the mathematical objects recognized at
that time.

To motivate our adherence to procedural issues, we note that there is
nothing wrong in principle with investigating \emph{pure ontology}.
However, practically speaking attempts by historians to gain insight
into Euler's ontology (as opposed to procedures) have a tendency, to
borrow Joseph Brodsky's comment in his introduction to Andrei
Platonov's novel \emph{The Foundation Pit}, to \emph{choke on their
own subjunctive mode}, as richly illustrated by an \emph{ontological}
passage that we quote in section \emph{Higher ontological order}
(Section~\ref{s34}).

The dichotomy of \emph{mathematical practice} versus \emph{ontology of
mathematical entities} has been discussed by a number of authors
including W.~Quine, who wrote: ``Arithmetic is, in this sense, all
there is to number: there is no saying absolutely what the numbers
are; there is only arithmetic.''  \cite[p.\;198]{Qu} 

For our purposes it will be more convenient to rely on Benacerraf's
framework.  \cite{Be65} pointed out that if observer~E learned that
the natural numbers ``are" the Zermelo ordinals
\[
\varnothing, \,  \{\varnothing\}, \,  \{\{\varnothing\}\}, \ldots,
\]
while observer~J learned that they are the von Neumann ordinals
\[
\varnothing, \, \{\varnothing\}, \, \{ \varnothing, \{\varnothing\}\},
\ldots
\]
then, strictly speaking, they are dealing with different things.
Nevertheless, observer~E's actual mathematical practice (and the
mathematical structures he is interested in) is practically the same
as observer~J's.  Hence, different ontologies may underwrite one and
the same practice.

For observer E, the entity~$0$ is \emph{not} an element of the
entity~$2$, while for observer J it is.  But for both of them the
\emph{relation}~$0<2$ holds.  Benacerraf's point is that although
mathematicians carry on their reasoning in terms of some objects or
others, the particular objects are not so important as the
\emph{relations} among those objects.  The relations may be the same,
even though the objects are different.

We would extend this insight beyond differences in set-theoretic
foundations and argue that even though Euler reasons about quantities
and Robinson reasons about sets (or types), they both agree, for
example, that~$a+dx=a$ for infinitesimal~$dx$ in a suitable
generalized sense of equality.  This is made precise in a hyperreal
framework via the standard part principle; see Section~\ref{21} for
more details.

This distinction relativizes the import of ontology in understanding
mathematical \emph{practice}.  A year after the publication of
Benacerraf's text, a related distinction was made by Robinson in
syntactic/procedural terms:
\begin{quote}
\ldots the theory of this book \ldots{} is presented, naturally,
within the framework of contemporary Mathematics, and thus appears to
affirm the existence of all sorts of infinitary entities.  However,
from a formalist point of view we may look at our theory syntactically
and may consider that what we have done is to introduce \emph{new
deductive procedures} rather than new mathematical entities.
\cite[p.\;282]{Ro66} (emphasis in the original)
\end{quote}

In short, we have, on the one hand, the \emph{ontological} issue of
giving a foundational account for the \emph{entities}, such as
infinitesimals and infinite integers, that classical infinitesimalists
may be working with.  On the other hand, we have their
\emph{procedures}, or \emph{inferential moves}, termed
\emph{syntactic} by Robinson.  What interests Euler scholars like
Laugwitz is not Euler's ontology but the syntactic \emph{procedures}
of his mathematical practice.  The contention that B-track
formalisations (see Section~\ref{12b}) provide better \emph{proxies}
for Euler's procedures and inferential moves than A-track
formalisations, is a \emph{methodological} or \emph{instrumentalist}
rather than an ontological or foundational matter.

To quote C. Pulte: ``Philosophy of science today should offer a more
accurate analysis to history of science without giving up its task --
not always appreciated by historians -- to uncover the basic concepts
and \emph{methods} which seem relevant for the understanding of
science in question.'' \cite[p.\;184]{Pu12} (emphasis added) Lagrange's
approach in his 1788 \emph{M\'echanique analytique} was remarkably
modern in its instrumentalism:
\begin{quote}
Neither are the metaphysical premises of his mechanics made explicit,
nor is there any epistemological justification given for the presumed
infallible character of the basic principles of mechanics. This is in
striking contrast not only to 17th-century foundations of mechanics
such as that of Descartes, Leibniz, and Newton but also to the
approaches of Lagrange's immediate predecessors, Euler, Maupertuis, or
d'Alembert \ldots{} 
In short, a century after Newton's \emph{Principia}, Lagrange's
textbook can be seen as an attempt to update the mathematical
principles of natural philosophy while abandoning the traditional
subjects of \emph{philosophia naturalis}.  In this \emph{special}
sense, the \emph{M\'echanique analitique} [sic] is also a striking
example of mathematical instrumentalism.  \cite[p.\;158]{Pu98}
(emphasis in the original)
\end{quote}
Two and a quarter centuries after Lagrange's instrumentalist approach,
perhaps a case can be made in favor of a historiography focusing on
methodological issues accompanied by an instrumentalist caution
concerning metaphysics and/or ontology of mathematical entities like
numbers and quantities.  This is in line with Pulte's insightful
comment made in the context of the study of rational mechanics in the
18th century: ``Euclideanism continues to be the ideal of science, but
it becomes a syntactical rather than a semantical concept of
science.'' \cite[p.\;192]{Pu12}

\subsection{A-track and B-track from Klein to Bos}
\label{12b}

The sentiment that there have been historically at least two possible
approaches to the foundations of analysis, involving dual methodology,
has been expressed by a number of authors.

In 1908, Felix Klein described a rivalry of two types of continua in
the following terms.  Having outlined the developments in real
analysis associated with Weierstrass and his followers, Klein pointed
out that ``The scientific mathematics of today is built upon the
series of developments which we have been outlining.  But an
essentially different conception of infinitesimal calculus has been
running parallel with this [conception] through the centuries''
\cite[p.\;214]{Kl08}.  Such a different conception, according to Klein,
``harks back to old metaphysical speculations concerning the structure
of the continuum according to which this was made up of \ldots{}
infinitely small parts'' (ibid.).  Thus according to Klein there is
not one but two separate tracks for the development of analysis:
\begin{enumerate}
\item[(A)] the Weierstrassian approach (in the context of an
\emph{Archimedean} continuum); and
\item[(B)] the approach with indivisibles and/or infinitesimals (in
the context of what could be called a \emph{Bernoullian} continuum).
\end{enumerate}
For additional details on Klein see Section~\ref{53}.

A similar distinction can be found in Henk Bos's seminal 1974 study of
Leibnizian methodology.  Here Bos argued that \emph{distinct}
methodologies, based respectively on (Archimedean) exhaustion and on
infinitesimals, are found in the work of 17th and 18th century giants
like Leibniz and Euler:
\begin{quote}
Leibniz considered two different approaches to the foundations of the
calculus; one connected with the classical methods of proof by
``exhaustion", the other in connection with a law of continuity
\cite[section 4.2, p.\;55]{Bos}.
\end{quote}
The first approach mentioned by Bos relies on an ``exhaustion''
methodology in the context of an Archimedean continuum.  Exhaustion
methodology is based on proofs by \emph{reductio ad absurdum} and the
ancient theory of proportion, which, as is generally thought today, is
based on the Archimedean axiom.%
\footnote{We note, in the context of Leibniz's reference to
Archimedes, that there are other possible interpretations of the
exhaustion method of Archimedes.  The received interpretation,
developed in \cite{Di87}, is in terms of the limit concept of real
analysis.  However, \cite [pp.\;280-290] {Wa85} developed a different
interpretation in terms of approximation by infinite-sided polygons.
The ancient exhaustion method has two components:
\begin{enumerate}
\item
\emph{geometric construction}, consisting of approximation by some
\emph{simple} figure, e.g., a polygon or a line built of segments,
\item
\emph{justification} carried out in the theory of proportion as
developed in \emph{Elements} Book V.
\end{enumerate}
In the 17th century, mathematicians adopted the first component, and
developed alternative justifications.  The key feature is the method
of exhaustion is the logical structure of its proof, namely
\emph{reductio ad absurdum}, rather than the nature of the background
continuum.  The latter can be Bernoullian, as Wallis' interpretation
shows.}

One way of formulating the axiom is to require that every positive
number can be added to itself finitely many times to obtain a number
greater than one.  The adjective \emph{Archimedean} in this sense was
introduced by O.~Stolz in the 1880s (see Section~\ref{s39}).  We will
refer to this type of methodology as the A-methodology.

Concerning the second methodology Bos notes: ``According to Leibniz,
the use of infinitesimals belongs to this kind of argument''
\cite[p.\;57]{Bos}.  We will refer to it as the B-methodology, in an
allusion to Johann Bernoulli (whose work formed the basis for
\cite{lh}), who, having learned an infinitesimal methodology from
Leibniz, never wavered from~it.

The Leibnizian laws such as the law of continuity mentioned by Bos in
the passage cited above, as well as the transcendental law of
homogeneity mentioned in \cite[p.\;33]{Bos}, find close procedural
analogs in Euler's work, and indeed in Robinson's framework.  The
transcendental law of homogeneity is discussed in Section~\ref{21}
and the law of continuity in Section~\ref{31a}.

In 2004, Ferraro appeared to disagree with Bos's \emph{dual track}
assessment, and argued for what he termed a ``continuous leap''
between (A-track) limits and (B-track) infinitesimals in Euler's work;
see Section~\ref{BF}.

\subsection{Mancosu and Hacking}
\label{s32}

To support our contention that there exist two distinct viable tracks
for the development of analysis, we call attention to Mancosu's
critique of G\"odel's heuristic argument for the inevitability of the
Cantorian cardinalities as the only plausible theory of the infinite.
G\"odel's argued that
\begin{quote}
the number of objects belonging to some class does not change if,
leaving the objects the same, one changes in any way whatsoever their
properties or mutual relations (e.g. their colors or their
distribution in space). \cite[p.\;254]{Go90}
\end{quote}
Mancosu argues that recent theories on numerosities undermine
G\"odel's assumption.  These were developed in \cite{BD} as well as
\cite{DF} and elsewhere.%
\footnote{A technical comment on numerosities is in order.  A
numerosity is a finitely additive measure-like function defined on an
algebra of sets, which takes values in the positive half of a
non-Archimedean ordered ring.  A numerosity is \emph{elementary} if
and only if it assigns the value~$1$ to every singleton in the domain,
so that the numerosity of any finite set is then equal to its number
of elements.  Therefore any elementary numerosity can be viewed as a
generalization of the notion of finite quantity.  Numerosities are
sometimes useful in studies related to Lebesgue-like and similar
measures, where they help to ``individualize" classically infinite
measure values, associating them with concrete infinitely large
elements of a chosen non-Archimedean ordered ring or field.  As a
concept of infinite quantity, numerosities have totally different
properties, as well as a totally different field of applications, than
the Cantorian cardinals.}
Mancosu concludes that
\begin{quote} 
having a different way of counting infinite sets shows that while
G\"odel gives voice to one plausible intuition about how to generalize
`number' to infinite sets there are coherent
alternatives. \cite[p.\;638]{Ma09}.
\end{quote}
Inspired in part by \cite{Ma09}, Ian Hacking proposes a distinction
between the \emph{butterfly model} and the \emph{Latin model}, namely
the contrast between a model of a deterministic biological development
of animals like butterflies, as opposed to a model of a contingent
historical evolution of languages like Latin.  For a further
discussion of Hacking's views see Section~\ref{hacking} below.

\subsection{Present-day standards}
\label{14}
Bos's comment on Robinson cited at the beginning of
Section~\ref{s21} is not sufficiently sensitive to the dichotomy of
practice (or procedures) versus ontology (or foundational account for
the entities) as analyzed in Section~\ref{11b}.  Leibnizian
procedures exploiting infinitesimals find suitable proxies in the
procedures in the hyperreal framework; see \cite{Re13} for a related
discussion in the context of Euler.  The relevance of such hyperreal
proxies is in no way diminished by the fact that set-theoretic
\emph{foundations} of the latter (``proof of the existence of the
entities,'' as Bos put it) were obviously as unavailable in the 17th
century as set-theoretic foundations of the real numbers.

In the context of his discussion of ``present-day standards of
mathematical rigor'', Bos writes: 
\begin{quote}
it is understandable that for mathematicians who believe that these
present-day standards are final, nonstandard analysis answers
positively the question whether, after all, \emph{Leibniz was right}
\cite[p.\;82, item~7.3]{Bos}. (emphasis added)
\end{quote}
The context of the discussion makes it clear that Bos's criticism
targets Robinson.  If so, Bos's criticism suffers from a strawman
fallacy, for Robinson specifically wrote that he did not consider set
theory to be the foundation of mathematics.  Being a formalist,
Robinson did not subscribe to the view attributed to him by Bos that
``present-day standards are final.''  Robinson expressed his position
on the status of set theory as follows: ``an infinitary framework such
as set theory \ldots{} cannot be regarded as the ultimate foundation
for mathematics'' \cite[p.\;45]{Ro69}; see also \cite[p.\;281]{Ro66}.
Furthermore, contrary to Bos's claim, Robinson's goal should not be
seen as showing that ``Leibniz was right'' (see above).  Rather,
Robinson's goal was to provide hyperreal proxies for the inferential
procedures commonly found in Leibniz as well as Euler and Cauchy.
Leibniz's procedures, involving as they do infinitesimals and infinite
numbers, seem far less puzzling when compared to their B-track
hyperreal proxies than from the viewpoint of the traditional A-track
frameworks; see Section~\ref{12b}.

\subsection{Higher ontological order}
\label{s34}

We wish to emphasize that we do \emph{not} hold that it is only
possible to interpret Euler in terms of modern formalisations of his
\emph{procedures}.  Discussions of Eulerian \emph{ontology} could
potentially be fruitful.  Yet some of the existing literature in this
direction tends to fall short of a standard of complete lucidity.
Thus, M.~Panza quotes Euler to the effect that ``Just as from the
ideas of individuals the idea of species and genus are formed, so a
variable quantity is the genus in which are contained all determined
quantities,'' and proceeds to explicate this as follows:
\begin{quote}
As constant quantities are determined quantities, this is the same as
claiming that a variable quantity is the genus in which are contained
all constant quantities. A variable quantity is thus a sort of a
formal characterization of quantity as such. Its concept responds to a
need for generality, i.e. a need of studying the essential properties
of any object of a certain genus, the properties that this object has
insofar as it belongs to such a genus. But, according to Euler, this
study has to have its own objects. In order to identify these objects,
it is necessary to sever these essential properties from any other
property that characterizes any object falling under the same
genus. If the genus is that of quantities, one has thus to identify
some objects that are not specific (and, a fortiori, particular)
quantities and pertain thus to a \emph{higher ontological order} than
that to which specific quantities pertain. \cite[pp.\;8-9]{Pa07}
(emphasis added).
\end{quote}
We are somewhat confused by this passage which seems to be ontological
in nature.  Since ontology is not our primary concern here (see
Section~\ref{11b}), we will merely propose further investigation into
the nature of variable quantities.
%

\section{Our reading of Leibniz and Euler}
\label{32}

The book \emph{Introductio in Analysin Infinitorum} \cite{ai} contains
remarkable calculations carried out in a framework where the basic
algebraic operations are applied to infinitely small and infinitely
large quantities.  

\subsection{Exponential function}
\label{s31}

In Chapter 7 on exponentials and logarithms expressed through series,
we find a derivation of the power series for the exponential
function~$a^z$ starting from the formula
\begin{equation}
\label{11}
a^\omega=1+k\omega .
\end{equation}
Here~$\omega$ is infinitely small, while~$k$ is finite.  Euler
specifically describes the infinitesimal~$\omega$ as being nonzero;
see Section~\ref{61}.  Euler then raises equation~\eqref{11} to the
infinitely great power~$i=\frac{z}{\omega}$ for a finite~$z$ to give
\begin{equation}
\label{12}
a^z= a^{i\omega}=(1+k\omega)^i .
\end{equation}
He then expands the right hand side of~\eqref{12} into a power series
by means of the binomial formula.  In the chapters that follow, Euler
finds infinite product decompositions for transcendental functions
(see Section~\ref{34} below where we analyze his infinite product
formula for sine).  In this section, we argue that the underlying
principles of Euler's mathematics are closer to Leibniz's than is
generally recognized.


\subsection{Useful fictions}
\label{uf}
\label{sync1}

We argue in this subsection that Euler follows Leibniz both
ontologically and methodologically.  On the one hand, Euler embraces
infinities as well-founded fictions; on the other, he distinguishes
assignable quantities from inassignable quantities.

The nature of infinitesimal and infinitely large quantities is dealt
with in Chapter~3 of \emph{Institutiones Calculi Differentialis}
\cite{ed}.  We cite Blanton's English translation of the Latin
original:
\begin{quote}
[e]ven if someone denies that infinite numbers really exist in this
world, still in mathematical speculations there arise questions to
which answers cannot be given unless we admit an infinite number.%
\footnote{In the original Latin this reads as follows: ``Verum ut ad
propositum revertamus, etiamsi quis neget in mundo numerum infinitum
revera existere; tamen in speculationibus mathematicis saepissime
occurrunt questiones, ad quas, nisi numerus infinitus admittatur,
responderi non posset.''  Note that the Latin uses the subjunctive
\emph{neget} (rather than \emph{negat}), which is the mode used for a
``future less vivid" condition: not ``even if someone denies" but
rather ``even if someone were to deny.''}
\cite[\S\,82]{Eu55}
\end{quote}
Here Euler argues that infinite numbers are necessary ``in
mathematical speculations'' even if someone were to deny ``their
existence in this world''.  Does this passage indicate that Euler
countenances the possibility of \emph{denying} that ``infinite numbers
really exist in this world''?  His position can be fruitfully compared
with that of the scholars of the preceding generation.  Indeed, the
scholars of the preceding generation disagreed on the issue of the
existence of infinitesimal quantities.  Bernoulli, l'H\^opital, and
Varignon staunchly adhered to the existence of infinitesimals, while
Leibniz adopted a more nuanced stance.  Leibniz's correspondence
emphasized two aspects of infinitesimal and infinite quantities: they
are
\begin{enumerate}
\item
\emph{useful fictions} and
\item
\emph{inassignable quantities}.
\end{enumerate}

It is important to clarify the meaning of the Leibnizian term
\emph{fiction}.  Infinitesimals are to be understood as \emph{pure
fictions} rather than \emph{logical fictions}, as discussed in
Section~\ref{s1}; see \cite{KS2}, \cite{KS1}, and \cite{SK}.
Furthermore, the work \cite{Je15} shows that Leibniz's strategy for
paraphrasing B-methods in terms of A-methods has to presume the
correctness of an infinitesimal inference (more precisely, an
inference exploiting infinitesimals), namely identifying the tangent
to a curve.  In the case of conic sections this succeeds because the
tangents are already known from Apollonius.  But for general curves,
and in particular for transcendental curves treated by Leibniz,%
\footnote{Leibniz applies his method in his \emph{de Quadratura
Arithmetica} to find the quadrature of general cycloidal segments
\cite[p.\;251]{Ed79}.  Here also the calculation exploits the family
of tangent lines.}
non-Archimedean infinitesimals remain an irreducible part of the
Leibnizian framework, contrary to~\cite[Chapter~5]{Is}.  This argument
is developed in more detail in \cite{Ba16}.

Similarly to Leibniz, Euler exploited the dichotomy of assignable
\emph{vs} inassignable, and mentioned the definition of infinitesimals
as being smaller than every assignable quantity, as well as the
definition of infinite numbers as being greater than every assignable
quantity; see \cite[p.\;17, 19, 20]{GKK}.  Thus, Euler writes: ``if~$z$
becomes a quantity less than any assignable quantity, that is,
infinitely small, then it is necessary that the value of the fraction
$1/z$ becomes greater than any assignable quantity and hence
infinite.'' \cite[\S\,90]{Eu55}.

Euler's wording in \cite[\S\,82]{Eu55}, making the usefulness of
infinite numbers \emph{independent} of their ``existence in this
world,'' suggests that his position is closer to a Leibnizian view
that infinitesimals are useful (or well-founded) fictions.  Euler goes
on to note that
\begin{quote}
an infinitely small quantity is nothing but a vanishing quantity, and
so it is really~$= 0$.%
\footnote{In the original Latin this reads as follows: ``Sed quantitas
infinite parva nil aliud est nisi quantitas evanescens, ideoque revera
erit~$= 0$.''  Note that the equality sign~``$=$'' and the
digit~``$0$'' are both in the original.  While Euler writes ``revera
erit~$= 0$" in \S\,83, in the next \S\,84 the formulation is ``revera
esse cyphram.''}
\cite[\S\,83]{Eu55}
\end{quote}

Euler's term \emph{nihil} is usually translated as \emph{nothing} by
Blanton.  However, in \emph{Introductio}, \S\,114, Blanton translates
``tantum non nihilo sit aequalis" as ``just not equal to zero" where
it should be ``just not equal to nothing" (see Section~\ref{61}).
Granted, ``equal to nothing" would sound awkward, but Euler seems to
distinguish it from ``equal to zero".  It is tempting to conjecture
that \emph{nihil} might be equivalent to ``exactly equal to zero",
whereas \emph{cyphra} is the term for a quantity whose only possible
assignable value is zero, or ``shadow zero".  Meanwhile in
\emph{Institutiones} \S\,84, Euler writes ``duae quaevis cyphrae ita
inter se sunt aequales, ut earum differentia fit nihil.''  This can be
translated as follows: ``two zeros are equal to each other, so that
there is no difference between them.''  This phrase is part of a
larger sentence that reads as follows in translation: ``Although two
zeros are equal to each other, so that there is no difference between
them, nevertheless, since we have two ways to compare them, either
arithmetic or geometric, let us look at quotients of quantities to be
compared in order to see the difference.'' \cite[p.\;51]{Eu55}.

This could be interpreted as saying that two instances of
\emph{cyphra} could be equal arithmetically but not geometrically.
The distinction between \emph{cyphra} and \emph{nihil} could
potentially give a satisfactory account for the Eulerian hierarchy of
zeros.

\subsection{Law of homogeneity from Leibniz to Euler}
\label{21}

As analyzed in Section~\ref{uf}, Euler insists that the relation of
equality holds between any infinitesimal and zero.  Similarly, Leibniz
worked with a generalized relation of ``equality'' which was an
\emph{equality up to} a negligible term.  Leibniz codified this
relation in terms of his \emph{transcendental law of homogeneity}
(TLH), or \emph{lex homogeneorum transcendentalis} in the original
Latin \cite{Le10b}.  Leibniz had already referred to the law of
homogeneity in his first work on the calculus: ``quantitates
differentiales, quae solae supersunt, nempe $dx$,~$dy$, semper
reperiuntur extra nominatores et vincula, et unumquodque membrum
afficitur vel per~$dx$, vel per~$dy$, servata semper \emph{lege
homogeneorum} quoad has duas quantitates, quomodocunque implicatus sit
calculus'' \cite{Le84} (emphasis added).  This can be translated as
follows: ``the only remaining differential quantities,
namely~$dx$,~$dy$, are found always outside the numerators and roots,
and each member is acted on by either~$dx$, or by~$dy$, always with
the law of homogeneity maintained with regard to these two quantities,
in whatever manner the calculation may turn out.''

The TLH governs equations involving differentials.  Bos interprets it
as follows:
\begin{quote}
A quantity which is infinitely small with respect to another quantity
can be neglected if compared with that quantity.  Thus all terms in an
equation except those of the highest order of infinity, or the lowest
order of infinite smallness, can be discarded.  For instance,
\begin{equation}
\label{adeq}
a+dx=a
\end{equation}
\[
dx+ddy=dx
\]
etc.  The resulting equations satisfy this \dots{} requirement of
homogeneity \cite[p.\;33]{Bos}
\end{quote}
(here the expression~$ddx$ denotes a second-order differential
obtained as a second difference).  Thus, formulas like Euler's 
\begin{equation}
\label{fo31}
a+dx=a
\end{equation}
(where~$a$ ``is any finite quantity''; \cite[\S\,\S\,86,87]{Eu55})
belong in the Leibnizian tradition of drawing inferences in accordance
with the TLH and as reported by Bos in formula~\eqref{adeq} above.
The principle of cancellation of infinitesimals was, of course, the
very basis of the technique, as articulated for example in \cite{lh}
(see also Section~\ref{52b}).  However, it was also the target of
Berkeley's charge of a logical inconsistency \cite{Be}.  This can be
expressed in modern notation by the conjunction
$(dx\not=0)\wedge(dx=0)$.  But the Leibnizian framework does not
suffer from an inconsistency of type~$(dx\not=0)\wedge(dx=0)$ given
the more general relation of ``equality up to''; in other words,
the~$dx$ is not identical to zero but is merely \emph{discarded} at
the end of the calculation in accordance with the TLH; see further in
Section~\ref{101}.

\subsection{Relations (\emph{pl.})~of equality}
\label{22}

What Euler and Leibniz appear to have realized more clearly than their
contemporaries is that there is more than one relation falling under
the general heading of ``equality''.  Thus, to explain formulas
like~\eqref{fo31}, Euler elaborated two distinct ways, arithmetic and
geometric, of comparing quantities.  He described the two modalities
of comparison in the following terms:

\begin{quote}
Since we are going to show that an infinitely small quantity is really
zero [\emph{cyphra}], we must meet the objection of why we do not
always use the same symbol 0 for infinitely small quantities, rather
than some special ones\ldots{} [S]ince we have two ways to compare
them [a more precise translation would be ``there are two modalities
of comparison''], either \emph{arithmetic} or \emph{geometric}, let us
look at the quotients of quantities to be compared in order to see the
difference.  (\cite{Eu55} \S\,84)
\end{quote}
Furthermore,
\begin{quote}
If we accept the notation used in the analysis of the infinite,
then~$dx$ indicates a quantity that is infinitely small, so that both
$dx =0$ and~$a\,dx=0$, where~$a$ is any finite quantity.  Despite
this, the \emph{geometric} ratio~$a\,dx: dx$ is finite, namely~$a:1$.
For this reason, these two infinitely small quantities,~$dx$
and~$a\,dx$, both being equal to~$0$, cannot be confused when we
consider their ratio.  In a similar way, we will deal with infinitely
small quantities~$dx$ and~$dy$.  (ibid., \S\,86, p. 51-52) (emphasis
added)
\end{quote}
Having defined the two modalities of comparison of quantities,
arithmetic and geometric, Euler proceeds to clarify the difference
between them as follows:
\begin{quote}
Let~$a$ be a finite quantity and let~$dx$ be infinitely small.
%
%
The arithmetic ratio of equals is clear: Since~$ndx =0$, we have
\begin{equation}
\label{26}
a \pm ndx - a = 0.
\end{equation}
On the other hand, the geometric ratio is clearly of equals, since
\begin{equation}
\label{32b}
\frac{a \pm ndx}{a} = 1.
\end{equation}
\cite[\S\,87]{Eu55}.
\end{quote}
While Euler speaks of distinct modalities of comparison, he writes
them down symbolically in terms of two \emph{distinct} relations, both
denoted by the equality sign ``$=$''; namely, \eqref{26} and
\eqref{32b}.  Euler concludes as follows:
\begin{quote}
From this we obtain the well-known rule that \emph{the infinitely
small vanishes in comparison with the finite and hence can be
neglected [with respect to it].} \cite[\S\,87]{Eu55} (emphasis in the
original)
\end{quote}
The ``well-known rule'' is an allusion to l'H\^opital's \emph{Demande
ou Supposition} discussed in Section~\ref{52b}.

Note that in the Latin original, the italicized phrase reads
\emph{infinite parva prae finitis evanescant, atque adeo horum
respectu reiici queant}.  The words ``with respect to it''
(\emph{horum respectu}) do not appear in Blanton's translation.  We
restored them because of their importance for understanding Euler's
phrase.  The term \emph{evanescant} can mean either \emph{vanish} or
\emph{lapse}, but the term \emph{prae} makes it read literally as
``the infinitely small vanishes \emph{before} (or \emph{by the side
of}) the finite,'' implying that the infinitesimal disappears
\emph{because of} the finite, and only once it is \emph{compared to}
the finite.

To comment on Euler's phrase in more detail, a possible interpretation
is that any motion or activity involved in the term \emph{evanescant}
does not indicate that the infinitesimal quantity is a dynamic entity
that is (in and of itself) in a state of disappearing, but rather is a
\emph{static} entity that changes, or disappears, only ``with respect
to'' (\emph{horum respectu}) a finite entity.  To Euler, the
infinitesimal has a different status depending on what it is being
compared to.  The passage suggests that Euler's usage accords more
closely with reasoning exploiting \emph{static} infinitesimals than
with dynamic limit-type reasoning.

Euler proceeds to present the usual rules going back to Leibniz,
L'H\^opital, and the Bernoullis, such as
\begin{equation}
\label{fo32}
a\,dx^m + b\,dx^n=a\,dx^m
\end{equation}
provided~$m<n$ ``since~$dx^n$ vanishes compared with~$dx^m$'' (ibid.,
\S\,89), relying on his geometric comparison.  Euler introduces a
distinction between infinitesimals of different order, and directly
\emph{computes} a ratio~$\frac{dx\pm dx^2}{dx}$ of two particular
infinitesimals by means of the calculation
\begin{equation}
\label{29}
\frac{dx\pm dx^2}{dx}=1\pm dx=1,
\end{equation}
assigning the value~$1$ to it (ibid., \S\,88).  Note that rather than
proving that the expression is equal to~$1$ (such \emph{indirect}
proofs are a trademark of the~$\epsilon, \delta$ approach), Euler
\emph{directly} computes (what would today be formalized as the
\emph{standard part} of) the expression.%
\footnote{To give an elementary example, the determination of the
limit~$\lim_{x\to 0} \frac{x+x^2}{x}$ in the~$\epsilon, \delta$
approach would involve first \emph{guessing} the correct
answer,~$L=1$, by using informal reasoning with \emph{small}
quantities; and then formally choosing a suitable~$\delta$ for
every~$\epsilon$ in such a way that~$\frac{x+x^2}{x}$ turns out to be
within~$\epsilon$ of~$L$ if~$|x|<\delta$.}
Euler combines the informal and formal stages by discarding the
higher-order infinitesimal as in \eqref{32b} and \eqref{29}.  Such an
inferential move is formalized in modern infinitesimal analysis in
terms of the standard part function or \emph{shadow}; see
Section~\ref{f15}.
Euler concludes:
\begin{quote}
Although all of them [infinitely small quantities] are equal to 0,
still they must be carefully distinguished one from the other if we
are to pay attention to their mutual relationships, which has been
explained through a geometric ratio (ibid., \S\,89).
\end{quote}
Like Leibniz in his \emph{Symbolismus} \cite{Le10b}, Euler considers
more than one way of comparing quantities.  Euler's
formula~\eqref{32b} indicates that his geometric comparison is
procedurally identical with the Leibnizian TLH (see
Section~\ref{21}): namely, both Euler's geometric comparison and
Leibniz's TLH involve discarding higher-order terms in the context of
a generalized relation of equality, as in \eqref{32b}
and~\eqref{fo32}.  

Note that there were alternative theories around 1700, such as the one
was proposed by Nieuwentijt.  Nieuwentijt's system, unlike Leibniz's
system, possessed only first-order infinitesimals with square zero
\cite{Ni}; \cite{Ver}; \cite[chapter~6]{Ma96}.  It is clear that the
Eulerian hierarchy of orders of infinitesimals follows Leibniz's lead.
%
%

Euler's geometric comparison was dubbed ``the principle of
cancellation'' in \cite[p.\;47]{Fe04}; see Section~\ref{e36} for a more
detailed discussion of Euler's \emph{zero infinitesimals}.

\subsection{Infinite product formula for sine}
\label{34}

In Section~\ref{12b} we analyzed a pair of approaches to interpreting
the work of the pioneers of analysis, namely the A-track in the
context of an Archimedean continuum, and the B-track in the context of
a Bernoullian continuum (an infinitesimal-enriched continuum).  We
explore a B-track framework as a proxy for the Eulerian
\emph{procedures}; here we leave aside the \emph{ontological} or
foundational issues, as discussed in Section~\ref{11b}.  We will
analyze specific procedures and inferential moves in Euler's oeuvre
and argue that the essential use he makes of both infinitesimals and
infinite integers is accounted for more successfully in a B-track
framework.

The fruitfulness of Euler's approach based on infinitesimals can be
illustrated by some of the remarkable applications he obtained.  Thus,
Euler derived an infinite product decomposition for the sine and sinh
functions of the following form:
\begin{eqnarray}
\label{eu1}
\sinh x &=& x\left(1+\frac{x^2}{\pi^2}\right)
\left(1+\frac{x^2}{4\pi^2}\right) \left(1+\frac{x^2}{9\pi^2}\right)
\left(1+\frac{x^2}{16\pi^2}\right) \,\dots\,\quad\quad \\[1ex]
\label{eu2}
\sin x &=& x\left(1-\frac{x^2}{\pi^2}\right)
\left(1-\frac{x^2}{4\pi^2}\right) \left(1-\frac{x^2}{9\pi^2}\right)
\left(1-\frac{x^2}{16\pi^2}\right) \,\dots\,\quad\quad
\end{eqnarray}
(see \emph{Introductio} \cite[\S\,155--164]{ai}).  Here~\eqref{eu2}
generalizes an infinite product formula for~$\frac{\pi}{4}$ (or
$\frac{\pi}{2}$) due to Wallis; see \cite[Proposition~191]{Wa}.
Namely, Wallis obtained the following infinite product:
\[
\prod_{n=1} \left(\frac{2n}{2n-1} \cdot \frac{2n}{2n+1}\right) =
\frac{2}{1} \cdot \frac{2}{3} \cdot \frac{4}{3} \cdot \frac{4}{5}
\cdot \frac{6}{5} \cdot \frac{6}{7} \cdot \frac{8}{7} \cdot
\frac{8}{9} \cdots = \frac{\pi}{2}.
\]
Evaluating Euler's product decomposition~$\frac{\sin x}{x} =
\prod_{n=1} \left(1 - \frac{x^2}{n^2\pi^2}\right)$ at
$x=\frac{\pi}{2}$ one obtains~$\frac{2}{\pi} = \prod_{n=1} \left(1 -
\frac{1}{4n^2}\right)$ or~$\frac{\pi}{2} = \prod_{n=1}
\left(\frac{4n^2}{4n^2 - 1}\right)$.  It follows that~$\frac{\pi}{2}=
\prod_{n=1} \left(\frac{2n}{2n-1}\cdot\frac{2n}{2n+1}\right)$, in
other words $\frac{\pi}{2} = \frac{2}{1} \cdot \frac{2}{3} \cdot
\frac{4}{3} \cdot \frac{4}{5} \cdot \frac{6}{5} \cdot \frac{6}{7}
\cdots$.

Euler also summed the inverse square
series:~$1+\frac14+\frac19+\frac1{16}+\ldots=\frac{\pi^2}6$; this is
the so-called \emph{Basel problem}.  This identity results
from~\eqref{eu2} by comparing the coefficient of~$x^3$ of the two
sides and using the Maclaurin series for sine.  This is one of Euler's
four solutions to the Basel problem; see \cite[p.\;111]{Sa}.

A common feature of these formulas is that Euler's computations
involve not only infinitesimals but also infinitely large natural
numbers, which Euler sometimes treats as if they were ordinary natural
numbers.

Euler's proof of the product decompositions~\eqref{eu1} and
\eqref{eu2} rely on infinitesimalist procedures that find close
proxies in modern infinitesimal frameworks.  In Appendix~\ref{31} we
present a detailed analysis of Euler's proof.

\subsection{Law of continuity}
\label{31a}

Euler's working assumption is that infinite numbers satisfy the same
rules of arithmetic as ordinary numbers.  Thus, he applies the
binomial formula to the case of an infinite exponent~$i$ without any
further ado in \cite[\S\,115]{ai}; see formula~\eqref{12} above.  The
assumption was given the following expression in 1755:
\begin{quote}
The analysis of the infinite, which we begin to treat now, is nothing
but a special case of the method of differences, explained in the
first chapter, wherein the differences are \emph{infinitely small},
while previously the differences were assumed to be
\emph{finite}. \cite[\S\,114]{Eu55} (emphasis added)
\end{quote}
The significance of this passage was realized by Bos (who gives a
slightly different translation; see Section~\ref{BF}).  Euler's
assumption is consonant with the Leibnizian \emph{law of continuity}:
\begin{quote}
il se trouve que les r\`egles du fini r\'eussissent dans
l'infini\ldots{} et que vice versa les r\`egles de l'infini
r\'eussissent dans le fini'' \cite{Le02}
\end{quote}
though apparently Euler does not refer explicitly to the latter in
this particular sense.  Robinson wrote:
\begin{quote}
Leibniz did say \ldots{} that what succeeds for the finite numbers
succeeds also for the infinite numbers and vice versa, and this is
remarkably close to our \emph{transfer of statements} from~$\R$
to~${}^\ast\R$ and in the opposite direction.  \cite[p.\;266]{Ro66}.
\end{quote}
On the \emph{transfer principle} see Section~\ref{e35}.  Euler treats
infinite series as polynomials of a specific infinite order (see
Section~\ref{81} for a discussion of the difference between finite and
infinite sums in Euler).  In the context of a discussion of the
infinite product
\begin{equation}
\label{33}
\begin{aligned}
\left(1+\frac{x}{i}+\frac{x^2}{4\pi^2}\right)
&
\left(1+\frac{x}{i}+\frac{x^2}{16\pi^2}\right)
\\&
\left(1+\frac{x}{i}+\frac{x^2}{36\pi^2}\right)
\left(1+\frac{x}{i}+\frac{x^2}{64\pi^2}\right)\cdots,
\end{aligned}
\end{equation}
where~$i$ is an infinite integer, Euler notes that a summand given by
an infinitesimal fraction~$\frac{x}{i}$ occurs in each factor.  One
may be tempted therefore to discard it.  The reason such an
infinitesimal summand cannot be discarded according to Euler, is
because it affects infinitely many factors:
\begin{quote}
through the multiplication of all factors, which are~$\frac{1}{2} i$
in number [$i$ being an infinitely large integer], there is a produced
term~$\frac{x}{2}$, so that~$\frac{x}{i}$ cannot be omitted
\cite[\S\,156]{ai}.
\end{quote}

In more detail, when one has a single factor, one can typically
neglect the infinitesimal~$\frac{x}{i}$.  However, in this case one
has~$\frac{i}{2}$ factors, and the linear term in the product will be
the sum of the linear terms in each factor.  This is one of the Vieta
rules that still holds when~$i$ is infinite by the \emph{law of
continuity}.  Altogether there are~$\frac{i}{2}$ factors, each of
which contains a linear term~$\frac{x}{i}$.  Therefore altogether one
obtains a contribution of~$\frac{i}{2} \cdot \frac{x}{i} =
\frac{x}{2}$, which is appreciable (noninfinitesimal) and therefore
cannot be neglected.

Euler's comment in 1748 shows that he clearly realizes that the
infinitesimal~$\frac{x}{i}$ present in each of the factors
of~\eqref{33} cannot be discarded at will.  While in 1755, the
\emph{preliminary} status of the infinitesimal is officially ``zero'',
in actual calculations Euler does not rely on such preliminary
declarations, as noted by Bos (see Section~\ref{61b} and
Section~\ref{BF}).

Leibniz's differentials~$dx$ were infinitesimals, and while Leibniz
did also consider non-infinitesimal differentials, he always denoted
them by the symbol~$(d)x$ rather than~$dx$;%
\footnote{Note that \cite{Bos} used the
notation~$\underline{\text{d}}x$ for Leibniz's~$(d)x$.}
see Section~\ref{sync1} for a discussion of Leibnizian infinitesimals.
There does not seem to be a compelling reason to think that
Euler's~$dx$'s were not infinitesimals, either.  Ferraro appears to
acknowledge this point when he writes: ``Euler often simply treats
differentials and infinitesimals as the same thing (for instance, see
Euler [1755, 70])'' \cite[p.\;35, note~2]{Fe04}.  Indeed, the
formula~$\omega=dx$ appears in \cite[\S\,118]{ed}.

Note that Euler explicitly refers to the \emph{number} of factors in
his infinite product, expressed by a specific infinite integer.
Similarly, when he applies the binomial formula~$(a+b)^i$ with an
infinite exponent~$i$, there is an implied \emph{final} term, or
\emph{terminal} summand, such as~$b^i$, though it never appears
explicitly in the formulas (see Section~\ref{81}).  We will analyze
Euler's proof in detail in Appendix~\ref{31}.

\subsection{The original rule of l'H\^opital}
\label{s36}

Euler's use of l'H\^opital's rule needs to be understood in its
historical context.  Most calculus courses today present the so-called
l'H\^opital's rule in a setting purged of infinitesimals.  It is
important to set the record straight as to the nature of the original
rule as presented by l'H\^opital in his \emph{Analyse des Infiniment
Petits pour l'Intelligence des Lignes Courbes}.

Two points should be kept in mind here.  First, L'H\^opital did
\emph{not} formulate his rule in terms of accumulation points, limits,
epsilons, and deltas, but rather in terms of infinitesimals:
\begin{quote}
Cela pos\'e, si l'on imagine une appliqu\'ee~$bd$ \emph{infiniment
proche} de~$BD$, \& qui rencontre les lignes courbes~$ANB$, 
$COB$ aux points~$f,g$; l'on aura~$bd=\frac{AB\times bf}{bg}$,
laquelle$^*$ ne diff\`ere pas de~$BD$ \cite[p.\;145]{lh}. (emphasis
added)
\end{quote}
A note in the right margin at the level of the asterisk following the
word \emph{laquelle} reads ``\;$^{*}$Art.~1.''  The asterisk refers
the reader to the following item:
\begin{quote}
\hfil I. Demande ou Supposition.

\noindent \ldots{} On demande qu'on puisse prendre indiff\'eremment
l'une pour l'autre deux quantit\'es qui ne different entr'elles que
d'une quantit\'e \emph{infiniment petite} \cite[p.\;2]{lh}. (emphasis
added)
\end{quote}
Clearly, Euler relied on l'H\^opital's original version of the rule
rather than any modern paraphrase thereof.  The original version of
l'H\^opital's rule exploited infinitesimals.  It seems reasonable
therefore that if one were to seek to \emph{understand} Euler's
procedures in a modern framework, it would be preferable to do so in a
modern framework that features infinitesimals rather than in one that
doesn't.

Our second point is that Euler's procedures admit a B-track
intrepretation in terms of an infinitesimal value of~$z$, and a
relation
\begin{equation}
\label{52}
\lambda\approx\frac{1-x^z}{z}
\end{equation}
of being infinitely close, or Euler's \emph{geometric comparison}; see
Section~\ref{22}.  These concepts are, on the one hand, closer to
Euler's world, and, on the other, admit rigorous proxies in the
context of a modern B-continuum (such as the hyperreals), namely the
relation~$\lambda = \text{st}\left( \frac{1-x^z}{z} \right)$ involving
the standard part function ``st''.  Arguably, the B-track
formula~\eqref{52} is a better proxy for \emph{understanding} Euler's
infinitesimal argument than is Ferraro's A-track formula~\eqref{51a}.

\subsection{Euclid's quantity}
\label{s38}

The classical notion of quantity is Euclid's \megethos{} (magnitude).
The general term \emph{magnitude} covers line segments, triangles,
rectangles, squares, convex polygones, angles, arcs of circles and
solids.  A general theory of magnitude is developed in the
\emph{Elements}, Book V.  In fact, Book V is a masterpiece of
deductive development.  By formalizing its definitions (see below the
formalisation of Definition V.4) and the tacit assumptions behind its
proofs, one can reconstruct Book~V and its~$25$ propositions as an
axiomatic theory.  \cite{Bec} and \cite[pp.\;101--122]{BM} provide
detailed sources for the axioms below in the primary source (Euclid).
See also \cite[pp.~118--148]{Mue} which mostly follows Beckmann's
development.  \cite{Heib83} is the standard modern edition of
\emph{Elements}.

As a result, Euclid's magnitudes of the same kind (line segments being
of one kind, triangles being of another, etc.) can be formalized as an
ordered additive semigroup with a total order~$<$ characterized by the
following five axioms:
\vspace{1ex}
\begin{itemize}\itemsep 0.5mm\itemsep 0.5mm
\item [E1]~$(\forall{x,y})(\exists n \in \mathbb{N}) [nx> y],$
\item [E2]~$(\forall{x,y})(\exists z) [ x<y \Rightarrow x+z=y],$
\item [E3]~$(\forall{x,y,z})[x < y \Rightarrow x+z <y+z],$
\item [E4]~$(\forall {x}) (\forall n \in \mathbb{N})(\exists y)[x=ny],$
\item [E5]~$(\forall x,y,z)(\exists v) [x:y::z:v]$.
\end{itemize}
Here axiom E1 formalizes \emph{Elements}, Definition~V.4.  More
specifically, Euclid's definition reads:
\begin{quote}
Magnitudes [such as~$a,b$] are said to have a ratio with respect to
one another which, being multiplied [i.e.,~$na$] are capable of
exceeding one another [i.e.,~$na> b$].
\end{quote}
The definition can be formalized as follows:~$(\forall a,b)(\exists
n)(na>b)$.
%
%
This reading of Euclid~V.4 is a standard interpretation among
historians; see \cite[pp.\;31--34]{Bec}; \cite[p.\;139]{Mue};
\cite[Section~II.3]{De16}.

The early modern mathematics developed largely without reference to
the Archimedean axiom.  Some medieval editions of \emph{Elements}
simply omitted the definition V.4; more precisely, they give
\emph{Proportion is a similarity of ratios} instead of definition V.4
of our modern editions; see \cite[p.\;137]{Grant}.  The same applies
to C.~Clavius' \emph{Euclidis Elementorum}, one of the most popular
{17th} century edition of \emph{Elements}; see
\cite[p.\;529]{Clavius}.

We do not find any explicit reference to the Archimedean axiom in the
works of Stevin, Descartes, Newton (though there is a mention of
Euclid's axiom in Leibniz's letter to l'H\^opital
\cite[p.\;288]{Le95a}), nor in the works of Euler.  Even the classical
constructions of the real numbers provided in 1872 by Heine, Cantor
and Dedekind contain no explicit mention of the Archimedean axiom, as
it was recognized as such only in 1880s by Stolz; see
Section~\ref{s39}.  The Archimedean axiom follows from the continuity
axiom (Dedekind axiom) and is equivalent to both the absence of
infinitesimals and the cofinality of the integers within the reals
defined in those constructions.  It took time for mathematicians to
understand the precise relation between the continuity axiom and the
Archimedean axiom.  It was not until 1901 that H\"older proved that
the continuity axiom (more precisely, Dedekind axiom) implies the
Archimedean axiom; see \cite{Ho} and \cite{Ho96}.

\subsection{Stolz and Heiberg}
\label{s39}

Otto Stolz in \cite{Sto85} rediscovered the Archimedean axiom for
mathematicians, making it one of his axioms for magnitudes.  The
Archimedean axiom was studied earlier in \cite{Sto83}, while
\cite{Sto85} was a popular and widely read book.  Stolz coined the
term \emph{Archimedean axiom}.  As the source of this axiom he points
out Archimedes' treatises \emph{On the sphere and cylinder} and
\emph{The quadrature of the parabola}.  As regards Euclid, Stolz
refers to books X and XII.  He does not seem to have noticed that
definition 4 of book V is related to the axiom of Archimedes.  Johan
L. Heiberg in his comment on the Archimedean axiom (lemma) cites
Euclid's definition~V.4 and observes that ``these are the same axioms"
\cite[p.\;11]{Heib81}.
%
%
Possibly as a result of his comment Euclid's definition V.4 is called
the Archimedean axiom.

At the end of the 19th century, Euclid's theory of magnitude was
revived by \cite{Sto85}, \cite{We}, and \cite{Ho}.  These authors
developed axiomatic theories of magnitude.  For a modern account of
these theories see \cite{Bl13}.  Despite certain differences, they all
accept axioms E1--E4 of Section~\ref{s38} as a common characterisation
of magnitude.  Instead of E5, some authors tend to use the Dedekind
axiom of continuity, which implies~E5.  H\"older was the first one to
show that E1 follows from E2--E4 and the Dedekind continuity axiom.

Thus, while axiom E1 is a feature of the classic and modern notion of
magnitude, it is absent from Euler's characterisation of quantity.
Moreover, Euler is explicit about the existence of infinite
quantities; see Section~\ref{121b}.

In the Eulerian context, a magnitude, or quantity, is not (yet) a
number.  Euler's quantities are converted to numbers once one
specifies an arbitrary quantity as the \emph{unit}, or \emph{unity}.
In addition to a unity, Euler needs a notion of a ratio.  Euler's
definition is similar to Newton's:
\begin{quote}
the determination, or the measure of magnitude of all kinds, is
reduced to this: fix at pleasure upon any one known magnitude of the
same species with that which is to be determined, and consider it as
the \emph{measure} or \emph{unit}; then, determine the proportion
[ratio] of the proposed magnitude to this known measure.  This
proportion [ratio] is always expressed by numbers; so that a number is
nothing but the proportion [ratio] of one magnitude to another
arbitrarily assumed as the unit \cite[\S~4]{Eu71}.
\end{quote}
However, neither Newton nor Euler provide a definition of
\emph{ratio}.  The term \emph{proportion} corresponds to the term
\emph{Verh\"{a}ltnis} (ratio) in the German edition of Euler's
\emph{Algebra}, and to \emph{rapport} (ratio) in the French edition
\cite{Eu07}.

\subsection{Euler on infinite numbers and quantities}
\label{121b}

Euler is explicit about the existence of infinite (and therefore
non-Archimedean) quantities and numbers:
\begin{quote}
not only is it possible to give a quantity of this kind, to which
increments are added without limit, a certain character, and with due
care to introduce it into calculus, as we shall soon see at length,
but also there exist real cases, at least they can be conceived,  in
which an infinite number actually exists \cite[\S~75]{Eu55}.
\end{quote}
Euler's important qualification ``at least they can be conceived''
with regard to the existence of infinite numbers is consistent with
the Leibnizian idea of them as \emph{useful fictions}; see
Section~\ref{21}.

\section{Critique of Ferraro's approach}
\label{primordial}

Ferraro's recent text on Euler seeks to steer clear of certain
interpretive approaches to Euler: ``My point of view is different from
that of some recent papers, such as [McKinzie-Tuckey 1997] and
[Pourciau 2001]. In this writing the authors recast the early
procedures directly in terms of the modern foundation of analysis or
interpret the earlier results in terms of modern theory of
non-standard analysis and understand the results in the light of this
later context'' \cite[p.\;2]{Fe12}.  

Ferraro's 2012 piece has significant textual overlap with his article
from 2004.  Here Ferraro asserts that ``one can see in operation in
their writings a conception of mathematics which is quite
\emph{extraneous} to that of Euler \ldots{} the attempt to specify
Euler's notions by applying modern concepts is only possible if
elements are used which are essentially \emph{alien} to them, and thus
Eulerian mathematics is transformed into something wholly different''
\cite[pp.\;51-52]{Fe04}; cf.~\cite[p.\;2]{Fe12} (emphasis added).  In
2004 Ferraro included two articles by Laugwitz in the list of such
allegedly ``extraneous'' and ``alien'' approaches: the article
\cite{Lau89} in \emph{Archive for History of Exact Sciences}, as well
as \cite{Lau92b}.

Ferraro's comments here betray an insufficient sensitivity to the
distinction analyzed in Section~\ref{11b}, namely, isolating
methodological concerns from obvious problems of ontology as far as
Euler's infinitesimals are concerned.  Granted, modern set-theoretic
frameworks, customarily taken to be an \emph{ontological} account of
the foundations of mathematics, are alien to Euler's world.  But is
Laugwitz's approach to Euler's \emph{methodology} really
``extraneous'' or ``alien'' to Euler?  Interpretive approaches seek to
clarify Euler's mathematical \emph{procedures} through the lens of
modern formalisations.  In the passage cited above, Ferraro appears
initially to reject such approaches, whether they rely on
modern~$\epsilon,\delta$ interpretations \`a la Weierstrass, or on
infinitesimal interpretations \`a la Robinson.  Yet in 2004, Ferraro
writes:
\begin{quote}
I am not claiming that 18th-century mathematics should be investigated
without considering modern theories.  Modern concepts are essential
for \emph{understanding} 18th-century notions and why these led to
\emph{meaningful} results, even when certain procedures, puzzling from
the present views, were used.  \cite[p.\;52]{Fe04} (emphasis added)
%
%
\end{quote}
(a similar passage appears in \cite[p.\;2]{Fe12}).  Thus, in the end
Ferraro does need modern theories to ``understand'' (as he puts it)
Euler, even though such procedures are ``meaningless'' to the latter.
Ferraro's position needs to be clarified, since any modern attempt to
\emph{understand} Euler will necessarily \emph{interpret} him, as
well.  While rejecting Laugwitz's interpretive approach to Euler,
Ferraro does seek to \emph{understand}, and therefore
\emph{interpret}, Euler by modern means.  To pinpoint the difference
between Laugwitz's interpretive approach (rejected by Ferraro) and
Ferraro's own interpretive approach, let us examine a sample of
Ferraro's reading of Euler.

\subsection{From l'H\^opital and Euler to epsilon and delta}
\label{52b}

Ferraro deals with an infinitesimal calculation in
\cite[pp.\;11-12]{Eu30} where Euler sought the value of
\[
\frac{1-x^{\frac{g}{f+g}}_{\phantom{I}}}{g}
\]
for~$f = 1$ and~$g = 0$ by applying l'H\^opital's rule to
$\frac{1-x^z}{z}$.  Ferraro proceeds to present the problem ``from a
modern perspective'' by analyzing the function~$f(z)=\frac{1-x^z}{z}$
and its behavior near~$z=0$ in the passage already cited in
Section~\ref{s21}, featuring the formula
\begin{equation}
\label{51a}
\lambda=\lim_{z\to c}f(z)
\end{equation}
(as already noted,~$c$ must be replaced by~$0$).  Here the formula
label~\eqref{51a} is added for later reference.

On the face of it, Ferraro merely explains what it means to a modern
reader to extend a function by continuity at a point where the
function is undefined.  However, Ferraro's presentation of a modern
explanation, with its talk of accumulation points, limits, epsilons,
and deltas, is firmly grounded in an A-track interpretation of the
Eulerian calculation using l'H\^opital's rule.  But why should one
seek to ``understand'' Euler using A-methodology?

In Section~\ref{s36} we placed Euler's use of l'H\^opital's rule in
its historiacal context.  What Ferraro presents here is
an~$(\epsilon,\delta)$ \`a la Weierstrass formalisation of Euler's
procedure.  He goes on to point out that such an approach would be
``meaningless'' to Euler.  Nevertheless, Ferraro goes on to make the
following remarkable claim:
\begin{quote}
there is something in common between the Eulerian procedure and the
modern one based upon the notion of limit: evanescent quantities and
endlessly increasing quantities were based upon an intuitive and
primordial idea of two quantities approaching each other.  I refer to
this idea as ``protolimit" to avoid any possibility of a modern
interpretation \cite[p.\;46]{Fe04}.
\end{quote}
Thus according to Ferraro, there is ``something in common between the
Eulerian procedure and the modern one,'' after all.  Ferraro's
\emph{protolimit} is intended to be different from the (A-track)
limit.  But shouldn't we rather interpret Eulerian infinitesimals in
terms of, say, a \emph{protoshadow}?  The term \emph{shadow} is
sometimes employed to refer to the (B-track) \emph{standard part
function}, discussed in Section~\ref{f15}.

\subsection{Shadow}
\label{f15}

In any totally ordered field extension~$E$ of~$\R$, every finite
element~$x\in E$ is infinitely close to a suitable unique element,
namely its standard part~$x_0\in\R$.%
\footnote{Indeed, via the total order, the element~$x$ defines a
Dedekind cut on~$\R$.  By the usual procedure, the cut specifies a
real number $x_0\in\R\subseteq E$.  The number~$x_0$ is infinitely
close to~$x\in E$.  The subring~$E_f\subseteq E$ consisting of the
finite (i.e., \emph{limited}) elements of~$E$ therefore admits a
map~$\text{st}:E_f\to\R,\;x\mapsto{}x_0$, called the \emph{standard
part function}, or \emph{shadow}, whose role is to \emph{round off}
each finite (limited)~$x$ to the nearest real~$x_0$.}

Ferraro finds fault with the standard part function as a tool in
intepreting Euler's equality~$\frac{i-1}{i}=1$.  Ferraro writes that
the equality ``should not be intended as~$\frac{i-1}{i}\approx 1$''
\cite[p.\;49]{Fe04} and provides the following clarification in
footnote~36 on page 49: ``By~$a\approx{}b$, I mean that the
difference~$a - b$ is an infinitesimal hyperreal number.''  The
criticism recurs in \cite[p.\;10]{Fe12} where the standard part
function is mentioned explicitly.  However, this criticism only raises
an issue if one assumes that Euler's equalities were not approximate
but, rather, exact equalities.  Such an assumption may be too
simplistic a reading of Euler's stance on \emph{arithmetic and
geometric comparisons}; see Section~\ref{32}.  See also
Appendix~\ref{31}, Step 5 and formula~\eqref{39b} where Euler wrote
that the term~$x^2/i^2$ is negligible in each of the factors of~$e^x -
e^{-x}$ \emph{only because} the number of the said factors is small
compared to~$i^2$.

Ferraro wrote that ``from the modern perspective, the problem of
extending the function'' is interpreted in terms of accumulation
points, A-track limits, epsilons, and deltas.

But couldn't we perhaps surmise, instead, that ``from a modern
perspective, the problem of extending the function may involve
infinitesimals, the relation of being infinitely close, and standard
part''?

Ferraro's claim that Eulerian infinitesimals ``were symbols that
represented a primordial and intuitive idea of limit''
\cite[p.\;34]{Fe04}, with its exclusive focus on the limit concept in
its generic meaning, tends to blur the distinction between the rival
Weierstrassian and modern-infinitesimal methodologies (see
Section~\ref{12b}).  Eulerian infinitesimals are intrinsically not
Archimedean but rather follow the methodology of his teacher
Bernoulli, co-founder with Leibniz of what we refer to as the
B-track. A better methodological proxy for Eulerian infinitesimals
than Ferraro's ``primordial limit'' is provided by a modern B-track
approach to analysis, fundamentally different from Ferraro's A-track
(proto)limit.

Meanwhile, Laugwitz sought to formalize Euler's procedures in terms of
modern infinitesimal methodologies.  It emerges that, while Ferraro's
own A-track reading is deemed ``essential for understanding
eighteenth-century notions and why these led to meaningful results''
as claimed in \cite[p.\;2]{Fe12}, Laugwitz's infinitesimal
interpretation is rejected by Ferraro as being both ``extraneous'' and
``alien'' to Euler's mathematics.  In short, Laugwitz's interpretation
does not fit Ferraro's Procrustean A-track way of, as he put it,
``understanding'' Euler \cite[p.\;2]{Fe12}.

\subsection{B-track reading in Felix Klein}
\label{53}

Laugwitz's interpretation accords with Felix Klein's remarks on the
dual tracks for the development of analysis as found in
\cite[p.\;214]{Kl08}.  In 1908, Felix Klein described a rivalry of the
dual approaches as we saw in Section~\ref{12b}.  Klein went on to
formulate a criterion for what would qualify as a successful theory of
infinitesimals.  A similar criterion was formulated in
\cite[pp.\;116-117]{Fran}.  For a discussion of the Klein--Fraenkel
criterion see \cite[Section~6.1]{KKM}.  The criterion was formulated
in terms of the mean value theorem.  Klein concluded:
\begin{quote}
I will not say that progress in this direction is impossible, but it
is true that none of the investigators have achieved anything
positive \cite[p.\;219]{Kl08}.
\end{quote}
Thus, the B-track approach based on notions of infinitesimals is not
limited to ``the work of Fermat, Newton, Leibniz and many others in
the 17th and 18th centuries,'' as implied by Victor J. Katz in
\cite{Ka14}.  Rather, it was very much a current research topic in
Felix Klein's mind.  See \cite{Eh06} for detailed coverage of the work
on infinitesimals around~1900.

Of course, Klein in 1908 had no idea at all of Robinson's hyperreal
framework as first developed in \cite{Ro61}.  What Klein was referring
to is the \emph{procedural} issue of how analysis is to be presented,
rather than the \emph{ontological} issue of a specific realisation of
an infinitesimal-enriched field in the context of a traditional set
theory; see Section~\ref{11b}.

Finally, we note that A-track readings of Euler tend to be
\emph{external} to Euler's \emph{procedures}, whereas infinitesimal
readings are \emph{internal},%
\footnote{This use of the term \emph{internal} is not to be confused
with its technical meaning in the context of enlargements of
superstructures; see \cite{Go}.}
in the sense that it provides proxies for \emph{both} the procedures
and the results of the historical infinitesimal mathematics.  This is
possible because modern infinitesimal procedures incorporate both
infinitesimals and infinite numbers as do Eulerian procedures.
Meanwhile, the Weierstrassian approach tends to provide proxies for
the results but not the procedures, since both infinitesimals and
infinite numbers have been eliminated in this approach.
%

\subsection{Hidden lemmas and principle of cancellation}
\label{e36}

Laugwitz argued that Euler's derivation of the power series expansion
of~$a^x$ contains a \emph{hidden lemma}, to the effect that a certain
\emph{infinite sum} of infinitesimals is itself infinitesimal under
suitable conditions; see \cite[p.\;210]{Lau89}.  Namely, let~$i$ be
infinite.  Consider Euler's formula
\[
\left( 1 + \frac{kz}{i} \right)^i = 1 + kz + \frac{i-1}{2i} k^2z^2 +
\frac{(i-1)(i-2)}{2i\cdot 3i} k^3z^3 + \ldots,
\]
or alternatively
\begin{equation}
\label{61c}
\left(1+\frac{kz}{i}\right)^i=1+\frac{i}{i} kz + \frac{i(i-1)}{i\cdot
2i} k^2z^2 + \frac{i(i-1)(i-2)}{i\cdot 2i\cdot 3i} k^3z^3 + \ldots
\end{equation}
There are infinitely many summands on the right.
\cite[\S\,115--116]{ai} goes on to make the substitutions
\begin{equation}
\label{62a}
\frac{i-1}{i}=1,\; \frac{i-1}{2i}=\frac{1}{2},\;
\frac{i-2}{3i}=\frac{1}{3},\ldots
\end{equation}
which he justifies by invoking the fact that~$i$ is infinite.  The
result is the exponential series.  The effect of these changes is
cumulative, since the products involved contain an ever increasing
number of factors.  Thus, one needs to make the substitution
\[
(i-1)(i-2)\cdots(i-n)=i^n
\]
for each finite~$n$ in the righthand side, but there are still
infinitely many summands affected.  Each of these substitutions
entails an infinitesimal change but there are infinitely many
substitutions involved in evaluating~\eqref{61c}.

Ferraro takes issue with Laugwitz's contention in the following terms:
\begin{quote}
It is evident that Laugwitz's remark arises from the interpretation of
$\frac{i-1}{i}= 1$ as~$\frac{i-1}{i}\approx 1$.  This interpretation
contrasts with the Eulerian statement that \hbox{$a+dx=a$} is an exact
equality and not an approximate one \cite[p.\;49]{Fe04}.
\end{quote}
Ferraro goes on to assert that, contrary to Laugwitz's claim, Euler
\begin{quote}
did not see gaps in the proof of [the series expansion of~$a^x$], and
this was due to the fact that he understood~$\frac{i-1}{i}= 1$ as a
formal equality involving fictitious entities \cite[p.\;50]{Fe04}.
\end{quote}
Indeed, if~$\frac{i-1}{i}= 1$ were an exact equality along with the
other expressions in~\eqref{62a}, the evaluation of the righthand side
of \eqref{61c} to the exponential series would be immediate and free
of any gaps, as Ferraro contends.

Alas, Ferraro underestimates Euler's perceptiveness here.  Ferraro
does not explain how an invocation of ``a formal equality involving
fictitious entities'' deflects Laugwitz's contention that Euler's
proof contains a \emph{hidden lemma}.  Ferraro's insistence on the
``exact equality''~$\frac{i-1}{i}= 1$ suggests that the infinitesimal
``error'' in~$\frac{i-1}{i}= 1$, or
\begin{equation}
\label{121}
1-\frac{1}{i}=1,
\end{equation}
is to be understood as exactly zero.  Declaring the infinitesimal
``error''~$\frac{1}{i}$ to be exactly zero would obviate the need for
justifying the hidden lemma, since an infinite sum of zeros is still
zero, or at any rate so Ferraro appears to interpret Euler's argument.
We will return to Ferraro's ``fictitious entities'' in
Section~\ref{52c}.

\subsection{Two problems with Ferraro's reading}
\label{61b}

There are two problems with Ferraro's claim that Euler is invoking an
exact equality with no infinitesimal error.  First, Euler explicitly
writes otherwise (see Section~\ref{31a} on the issue of disappearing
infinitesimals), and in fact in the calculation under discussion,
Euler exploits the relation~$z=\omega i$ \cite[\S\,115]{ai} with
infinitesimal~$\omega$ and finite~$z$, which would be quite impossible
if~$\omega$ were literally zero.

The second problem is that, as Ferraro himself noted in his recent
text, Euler expressed the integral as ``the sum of an infinite number
of infinitesimals'' \cite[p.\;10]{Fe12}.
%
%
Euler expresses the integral in terms of the expression
\begin{equation}
\label{39}
\alpha (A+A'+A''+A'''+\ldots+X)
\end{equation}
in \cite[Chapter~VII, p.\;184]{Eu68}, where~$\alpha$ is an
infinitesimal step of a suitable partition,
while~$A$,~$A'$,~$A''$,~$A''',\ldots$ are the (finite) values of the
integrand at (infinitely many) partition points.  The quantities
\[
\alpha A,\;\alpha A',\;\alpha A'',\;\alpha A''',\ldots
\]
are still infinitesimal, and therefore would be exactly zero, so that
their infinite sum~\eqref{39} would be paradoxically zero as well.
Thus, such a reading of Euler's reasoning attributes to him an
alarming paradox not dealt with in Ferraro's approach.

Ferraro mentions Euler's interpretation of the integral as an infinite
sum of infinitesimals in \cite[p.\;50, footnote~39]{Fe04}, but fails to
explain how the paradox mentioned in the previous paragraph could be
resolved (other than implying that infinitesimals are sometimes zero,
and sometimes not).

In sum, we agree with Bos' evaluation of Euler's preliminary remarks
on ``infinitesimals as zeros'' as being at variance with his actual
mathematical practice (see Section~\ref{BF}).  It is unlikely that a
literal interpretation of Euler's preliminary remarks (that the
infinitesimal is exactly zero) could give a fruitful way of
interpreting Euler's mathematics.  Ferraro's rejection of Laugwitz's
analysis of Euler's argument in terms of a \emph{hidden lemma}
(requiring further justification) is therefore untenable.

\subsection{Generality of algebra}
\label{e35}

It was known already to Cauchy that some of Euler's doctrines are
unsatisfactory.  More specifically, Cauchy was critical of Euler's and
Lagrange's \emph{generality of algebra}, to the effect that certain
relations involving variable quantities are viewed as being valid even
though they can fail for certain specific values of the variables.  By
the time mathematicians started analyzing Fourier series in the 1820s
it became clear that some applications of the \emph{generality of
algebra} are untenable.  Cauchy specifically rejects this principle in
the introduction to his \emph{R\'esum\'e des Le\c cons} \cite{Ca23}.

In the context of a discussion of Euler's principle of the
\emph{generality of algebra}, Ferraro notes that the Eulerian
``general quantity''
\begin{quote}
was represented by graphic signs which were manipulated according to
appropriate rules, which were the same rules that governed geometrical
quantities or true numbers \cite[p.\;43]{Fe04}.
\end{quote}
The idea that ``the same rules'' should govern ideal/fictional numbers
and ``true numbers'' is consonant with the Leibnizian \emph{law of
continuity}.  The latter is arguably a fruitful methodological
principle.  It was formalized as the \emph{transfer principle} in
Robinson's framework.%
\footnote{The \emph{transfer principle} is a type of theorem that,
depending on the context, asserts that rules, laws or procedures valid
for a certain number system, still apply (i.e., are ``transfered'') to
an extended number system.  Thus, the familiar
extension~$\Q\subseteq\R$ preserves the properties of an ordered
field.  To give a negative example, the
extension~$\R\subseteq\R\cup\{\pm\infty\}$ of the real numbers to the
so-called \emph{extended reals} does not preserve the properties of an
ordered field.  The hyperreal extension $\R\subseteq\astr$ preserves
\emph{all} first-order properties, such as the identity~$\sin^2 x +
\cos^2 x =1$ (valid for all hyperreal~$x$, including infinitesimal and
infinite values of~$x\in\astr$).  For a more detailed discussion, see
the textbook \emph{Elementary Calculus} \cite{Ke86}.}

Meanwhile, Ferraro fails to distinguish between, on the one hand, a
historically fruitful \emph{law of continuity}, and, on the other, the
\emph{generality of algebra} that was found to be lacking in the 19th
century, as he continues in the next sentence:
\begin{quote}
The principle of the generality of algebra held: the rules were
applied in general, regardless of their conditions of validity and the
specific values of quantity (ibid.).
\end{quote}

Cauchy's critique of Euler's principle of the \emph{generality of
algebra} is well known to historians; it is an uncontroversial
statement that certain elements of Euler's oeuvre need to be
reinterpreted if one is to develop a consistent interpretation
thereof.  Another such element is the zero infinitesimal, as discussed
in Section~\ref{63}.

\subsection{Unsettling identity}
\label{63}

The claim that the infinitesimal is exactly \emph{equal} to zero
occasionally does appear in Euler's writing, such as in the
\emph{Institutiones} in reference to~$dx$.  On the other hand, Euler
specifically discusses varieties of the notion of \emph{equality},
with the \emph{geometric} notion being similar to the generalized
relation of equality implied in the Leibnizian \emph{transcendental
law of homogeneity} (see Section~\ref{22}).  Even though at times
Euler insists that his equality is \emph{exact equality}, at other
times he does envision more general modalities of comparison.  Ferraro
himself implicitly acknowledges this when he describes Euler's
equality~$a+dx=a$ as a ``principle of cancellation''
\cite[p.\;47]{Fe04}.  The term \emph{principle of cancellation} would
appear to imply that there is something to \emph{cancel}: not merely
an \emph{exact zero}, but a nonzero infinitesimal~$dx$.

On an even more basic level, if for infinite~$i$ one
has~$\frac{1}{i}=0$ as in~\eqref{121}, then multiplying out by~$i$ we
obtain~$1=0 \times i$, but~$0$ times any number is still~$0$, so that
we would obtain an unsettling identity~$1=0$ (at least if we interpret
``$=$'' as literal equality), in addition to the paradox with the
integral mentioned above.

Similarly, Euler seeks to divide by an infinitesimal~$dx$ so as to
obtain the differential ratio~$\frac{dy}{dx}$.  It follows that~$dx$
cannot be an exact zero if one is to have any hope for a consistent
account of Euler's procedures.  A notion of \emph{zero infinitesimal}
interpreted literally is arguably as problematic as some aspects of
the principle of the \emph{generality of algebra} already found to be
lacking by Cauchy (see Section~\ref{e35}).  It can be reinterpreted in
terms of the distinction between \emph{cyphra} and \emph{nihil} as
discussed in Section~\ref{uf}.

\subsection{Fictitious entities}
\label{52c}
\label{61}

Ferraro described the Eulerian substitution
$\frac{i-1}{i}= 1$ as a ``formal equality involving fictitious
entities'' in his text \cite[p.\;50]{Fe04}.  It is not entirely clear
how such an evocation of fictionality resolves the delicate
mathematical problem posed by this substitution.

Two scientific generations earlier, Leibniz described infinitesimals
as ``useful fictions,'' yet he did not think that just because
infinitesimals are ``fictions'' one is allowed to set them equal to
zero at will.

Ferraro further claimed that ``The use of fictions made Eulerian
mathematics extremely different from modern mathematics" in his text
\cite[p.\;64]{Fe07}.  The plausibility of the claim depends on
equivocation on procedure/ontology as discussed in Section~\ref{11b}.
The fact that what Ferraro has in mind here is ontology is made clear
on the previous page where he writes: ``The rules that Euler uses upon
[\emph{sic}] infinite and infinitesimal quantities constitute an
immediate extrapolation of the behaviour of a finite variable~$i$
tending to~$\infty$ or~$0$'', \cite[p.\;63]{Fe07} and adds: ``what is
wholly missing is the \emph{complex construction} of
${}^\ast\mathbb{R}$ and the assumptions upon which it is based.''
(ibid.) (emphasis added) But this ontological complaint is utterly
irrelevant to procedure.

The infinitesimals in Leibniz and Euler may have been \emph{fictions}.
However, they were not \emph{fictive} or purely rhetorical, as Ferraro
appears to imply.  Rather, they pose subtle issues of interpretation
that are not resolved by an appeal to ``formal equality involving
fictitious entities;'' see further in Section~\ref{eight}.

On occasion, Euler specifically wrote that his infinitesimals are
\emph{unequal} to zero: ``Let~$\omega$ be a number infinitely small,
or a fraction so tiny that it is just not equal to zero (\emph{tantum
non nihilo sit aequalis})'' \cite[\S\,114]{aie}.  This passage refers
specifically to the infinitesimal~$\omega$ in formula~\eqref{11} used
in the derivation of the power series of the exponential function (see
Section~\ref{s31} on the exponential function), showing that error
estimates are indeed required even if one takes Euler literally.

\subsection{Finite, infinite, and hyperfinite sums}
\label{81}

Euler's use of infinite integers and their associated infinite
products (such as the product decomposition of the sine function) were
interpreted in Robinson's framework in terms of hyperfinite
expressions.  Thus, Euler's product of~$i$-infinitely many
%
%
factors in~\eqref{eu2} is interpreted as a hyperfinite product in
\cite [formula~(9), p.\;74] {kr}.  A hyperfinite formalisation of
Euler's argument involving infinite integers and their associated
products illustrates the successful formalisation of the
\emph{arguments} (and not merely the \emph{results}) of classical
infinitesimal mathematics.
%

In a footnote on 18th century notation, Ferraro presents a novel claim
that ``for 18th-century mathematicians, there was no difference
between finite and infinite sums'' \cite[footnote~8, p.\;294]{Fe98}.
Far from being a side comment, the claim is emphasized a decade later
in the Preface to his book: ``a distinction between finite and
infinite sums was lacking, and this gave rise to formal procedures
consisting of the infinite extension of finite procedures''
\cite[p.\;viii]{Fe08}.  The clue to decoding Ferraro's claim is found
in the same footnote, where Ferraro distinguishes between sums
featuring a final term after the ellipsis, such as
\begin{equation}
\label{41}
a_1+a_2+\ldots+a_n,
\end{equation}
and ``infinite sums'' without such a final term, as in
\begin{equation}
\label{42}
a_1+a_2+\ldots+a_n+\ldots
\end{equation}

Note that A-track syntax is unable to account for \emph{terminating}
infinite expressions which routinely occur in Euler.  To be sure,
Euler does \emph{not} use his infinite~$i$ as a final index in
infinite sums of type~$a_1+\ldots+a_i$ common in modern infinitesimal
frameworks.  However, his binomial expansions with exponent~$i$ play
the same role as the modern infinite sums~$a_1+\ldots+a_i$.  The final
term~$a_i$ is hinted upon by means of Euler's notation ``\&c.'' but
does not appear explicitly.  Nonetheless, procedurally speaking his
infinite sums play the same role as the modern~$a_1+\ldots+a_i$.

From an A-track viewpoint, a terminating sum~\eqref{41} is necessarily
a finite one, whereas only expressions of the form~\eqref{42} ending
with an ellispis allow for a possibility of an ``infinite sum.''  No
other option is available in the A-track; yet Euler appears recklessly
to write down \emph{infinite} terminating expressions, as in the proof
of the product formula for sine (see Section~\ref{31a} for a
discussion of terminal summands in infinite sums in Euler).

Meanwhile, the B-track approach allows one to account both for Euler's
\emph{infinite} integer~$i$ and for \emph{terminating} expressions
containing~$i$ terms (see Appendix~\ref{31} for an instance of Euler's
use of polynomials of infinite degree).  Euler discusses the
difference between finite and infinite sums in
\emph{Introductio}~\cite[\S\,59]{ai}.  Terminating \emph{infinite}
sums are easily formalized in Robinson's framework in terms of
hyperfinite expressions (see Appendix~\ref{31}).

In a subsequent article, Ferraro and Panza write: ``Power series were
conceived of as quasi-polynomial entities (that is, mere
\emph{infinitary} extensions of polynomials)'' \cite[p.\;20]{FP})
(emphasis added), but don't mention the fact that such an extension
can be formalized in terms of hyperfinite expressions, perhaps out of
concern that this may be deemed ``alien'' or ``extraneous'' to Euler.

Euler's formula~$a^{i\omega}=(1+k\omega)^i$ is analyzed in
\cite[p.\;48]{Fe04} (a similar analysis appears in \cite[p.\;9]{Fe12} in
a section entitled ``Non-standard analysis and Eulerian
infinitesimals'').  Ferraro reformulates Euler's formula in terms of
modern \emph{Sigma} notation as follows:
\[
a^{i\omega}=(1+k\omega)^i = \sum_{r=0}^\infty {i \choose r}
(k\omega)^r .
\]
The formula
\begin{equation}
\label{43}
a^x=\sum_{r=0}^{\infty}\frac{1}{r!} (kx)^r
\end{equation}
appears in \cite[Formula~(2), p.\;48]{Fe04} and is attributed to
Euler.  The Sigma notation~$\sum_{r=0}^\infty$ appears several times
in Ferraro's analysis and is clearly not a misprint; it appears again
in \cite[p.\;48, 54]{Fe07}.  Similarly, Ferraro exploits the modern
notation~$\sum_{i=1}^\infty a_i$ for the sum of the series, in
\cite[p.\;5]{Fe08}, while discussing late 16th (!) century texts of
Vi\`ete.  The Sigma notation~\eqref{43} is familiar modern notation
for infinite sums defined via the modern concept of limit in a
Weierstrassian context.  Note that formula~\eqref{43} attributed by
Ferraro to Euler involves assigning a sum to the series, namely~$a^x$,
and therefore is not merely a formal power series.  The summation of
an infinite series via the concept of limit (namely, limit of the
sequence of partial sums) is not accessory but rather a \emph{sine qua
non} aspect of such summation (alternatively, one could take the
standard part of a hyperfinite sum, but such an approach is apparently
not pursued by Ferraro).  The symbol~$\infty$ appears in~\eqref{43} as
a kind of subjunctive.  It has no meaning other than a reminder that a
limit was taken in the definition of the series.  In modern notation,
the symbol~$\infty$ does not stand for an infinite integer (contrary
to the original use of this symbol by Wallis).

Thus, Ferraro reformulates Euler's calculation using the Sigma
notation for infinite sum, including the modern somewhat subjunctive
use of the superscript~$\infty$.  However, such a procedure is
extraneous to Euler's mathematics, since Euler specifically denotes
the (infinite) power by~$i$.  Applying the binomial formula with
exponent~$i$, one would obtain, not Ferraro's~\eqref{43}, but rather a
sum from~$0$ to~$i$, namely
\begin{equation}
\label{44}
\sum_{r=0}^i {i \choose r} (kx)^r .
\end{equation}
Euler's derivation of the exponential series is analyzed in more
detail in Section~\ref{e36}.  Infinite sums of type~\eqref{44} are
perfectly meaningful when interpreted in Robinson's framework (see
Appendix~\ref{31}).  Ferraro's anachronistic rewriting of Euler's
formula betrays a lack of sensitivity to the actual mathematical
content of Euler's work.

\subsection{Bos--Ferraro differences}
\label{BF}

In this section we compare Ferraro's take on Euler with the approaches
by other scholars, more compatible with our reading of Euler.  We will
first compare the approaches of Bos and Ferraro to Euler scholarship,
and then those of Ferraro and Laugwitz.  Bos summarized Euler's
preliminary discussion of infinitesimals in the following terms:
\begin{quote}
Euler claimed that infinitely small quantities are equal to zero, but
that two quantities, both equal to zero, can have a determined ratio.
This ratio of zeros was the real subject-matter of the differential
calculus \cite[p.\;66]{Bos}.
\end{quote}
Bos goes on to note that Euler's preliminary discussion is at variance
with Euler's actual mathematical practices even in the
\emph{Institutiones} (and not merely in the \emph{Introductio} as
discussed in Section~\ref{31a}), where the properties of the
infinitely small are similar to those of finite differences:
\begin{quote}
After having treated, in the first two chapters, the theory of finite
difference sequences, he defined the differential calculus as the
calculus of infinitesimal differences:
\begin{quote}
The analysis of infinites, with which I am dealing now, will be
nothing else than a special case of the method of differences
expounded in the first chapter, which occurs, when the differences,
which previously were supposed finite, are taken infinitely small
\cite[\S 114]{ed},
\end{quote}
which is rather at variance with his remarks quoted above, a
\emph{contradiction} which shows that his arguments about the
infinitely small did not really influence his presentation of the
calculus \cite[pp.\;67--68]{Bos} (emphasis added)
\end{quote}

Before analyzing Ferraro's reaction to this position, we note that
Bos' focus on Euler's ``presentation of the calculus'' indicates a
concern for methodological issues related to the nature of Euler's
\emph{procedures}, rather than focusing on the \emph{ontological}
nature of the objects (the infinitely small) that Euler utilizes, in
line with the distinction between procedure and ontology that we
explored in Section~\ref{11b}.

Ferraro disagrees with Bos' perception of a ``contradiction'' in
Euler's writing:
\begin{quote}
According to Bos, there is ``a contradiction which shows that his
arguments about the infinitely small did not really influence his
presentation of calculus" [Bos, 1974, 68-69].  However, I would argue
that one may see a contradiction in the \emph{Institutiones} only if,
in contrast to Euler, [1] one distinguishes between \emph{limits} and
\emph{infinitesimals} and [2] neglects the nature of evanescent
quantities as fictions, [3] the role of formal manipulations and
[4]~the absence of a separation between semantics and syntax in the
Eulerian calculus \cite[p.\;54]{Fe04} (emphasis and numerals [1], [2],
[3], [4] added).
\end{quote}
Ferraro appears to suggest that Bos' position is problematic with
regard to the four items enumerated above.  We will not analyze all
four, but note merely that in his item
\begin{equation}
\label{71b}
\text{[1] ``one distinguishes between limits and infinitesimals''},
\end{equation}
Ferraro commits himself explicitly to the position that
``distinguish[ing] between limits and infinitesimals,'' as Bos does,
is an inappropriate approach to interpreting Euler.  Rather, Ferraro
sees a conceptual continuity between limits and infinitesimals in
Euler, or more precisely what he refers to as a ``continuous leap''
(see Section~\ref{71}).

We argue that Bos's position on this aspect of Euler's \emph{oeuvre}
is more convincing than Ferraro's.  Note that Euler's insistence on
the similarity of the properties of the finite and infinitesimal
differences, in the passage cited by Bos, is consonant with a
Leibnizian \emph{law of continuity}, which requires two types of
quantities to be compared: assignable and inassignable (e.g.,
infinitesimal); see Section~\ref{21}.

\subsection{Was Euler ambiguous or confused?}
\label{71}

Ferraro postulated a conceptual continuity between limits and
infinitesimals in Euler's work, as expressed in Ferraro's
comment~\eqref{71b} meant to be critical of Bos' position.


Ferraro's criticism of Bos' approach emanates from Ferraro's tendency
to blur the distinction between A-track and B-track approaches.  A
further attempt to blur this distinction is found in Ferraro's
``continuous leap'' comment:
\begin{quote}
Eulerian infinitesimals \ldots when interpreted using the conceptual
instruments available to modern mathematics, seem to be an
\emph{ambiguous} mixture of different elements, a \emph{continuous
leap} from a vague idea of \emph{limit} to a \emph{confused} notion of
\emph{infinitesimal}'' \cite[p.\;59]{Fe04} (emphasis added)
\end{quote}
Ferraro's comment appears in the ``Conclusion'' section in 2004.  A
virtually identical comment appears in the abstract in 2012, and yet
again in the ``Conclusion'' section in \cite[p.\;24]{Fe12}.

%
%

We argue however that Euler was far less ``ambiguous'' or ``confused''
than is often thought.  Ferraro claims that when we allow our
interpretation of Euler to be informed by modern mathematical
concepts, we have no choice but to see Euler as fluidly moving from
vague limits to confused infinitesimals.  Let us now compare the
interpretations by Ferraro and Laugwitz.

Ferraro's opposition to Laugwitz's interpretation is based on a
conflation of ontology and practice (see Section~\ref{11b}).  Laugwitz
is not trying to read ontological foundations based on modern theories
into Euler (which would indeed be ``alien" to Euler's notions, to
borrow Ferraro's terminology), but is rather focusing on Euler's
mathematical practice.

Furthermore, Ferraro's own reading, with its emphasis on alleged
continuity between limits and infinitesimals, is not sufficiently
sensitive to the distinction between dual approaches (as analyzed by
both Klein and Bos), which we refer to as A-track and B-track
approaches.

Laugwitz's interpretation showed that drawing upon modern concepts
allows us to see Euler's reasoning as clear and incisive.  Indeed, we
know since \cite{Ro61} that Felix Klein's hunch concerning the
\emph{dual} approaches to the foundations of analysis in
\cite[p.\;214]{Kl08} was right on target (see Section~\ref{53}).  In
short, Ferraro assimilates two distinct approaches to the problem of
the continuum without historical or mathematical evidence.

\subsection{Rhetoric and modern interpretations}
\label{eight}

In his 2004 article \cite[p.\;51, footnote~46]{Fe04} Ferraro sought to
enlist the support of \cite[Appendix~2]{Bos} for his (Ferraro's)
opposition to interpretating Euler in terms of modern theories of
infinitesimals.  However, Henk Bos himself has recently distanced
himself from the said Appendix~2 (part of his Doctoral thesis) in a
letter sent in response to a question from one of the authors of the
present text:
\begin{quote}
An interesting question, what made me reject a claim some 35 years
ago?  I reread the appendix and was surprised about the self assurance
of my younger self.  I'm less definite in my opinions today -- or so I
think.  You're right that the appendix was not sympathetic to
Robinson's view.  Am I now more sympathetic? If you talk about
``historical continuity" I have little problem to agree with you,
given the fact that one can interpret continuity in historical
devlopments in many ways; even revolutions can come to be seen as
continuous developments \cite{Bos12}.
\end{quote}
The letter is reproduced with the author's permission.  The
shortcomings of Bos's Appendix~2 are analyzed in detail in
\cite[Section~11.3]{KS1} and in Section~\ref{14} here.  The
clarification provided by Bos in 2010 weakens the claim of Bos's
support for Ferraro's position on Robinson.  Ferraro claims that
\begin{quote}
[Laugwitz and other] commentators use notions such as set, real
numbers, continuum as a set of numbers or points, functions as
pointwise relations between numbers, axiomatic method, which are
modern, not Eulerian \cite[p.\;51]{Fe04}.
\end{quote}
Certainly, sets, real numbers, the punctiform continuum, and the
modern notion of function are not Eulerian concepts.  But has Laugwitz
really committed the misdemeanors attributed to him by Ferraro?
Ferraro does not provide any evidence for his claim, and one searches
in vain the two articles \cite{Lau89} and \cite{Lau92b} cited by
Ferraro for clues of such misdemeanors.  On the contrary, Laugwitz
warns the reader: ``But one should have in mind that such concepts did
not appear before set theory was established'' \cite[p.\;242]{Lau89};
and again:
\begin{quote}
Modern mathematicians should find of interest the fact that he
[Cauchy] succeeded by using only very few concepts of an intensional
quality, whereas we have become accustomed to using a great many
extensional concepts based on set theory (ibid.).
\end{quote}
Laugwitz is clearly aware of the point that modern set theory is
\emph{alien} to Euler's ontology.  However, as discussed in
Section~\ref{11b}, Laugwitz is concerned with Euler's
\emph{procedures} rather than his ontology.  Ferraro has surely
committed a strawman fallacy in describing Laugwitz's scholarship as
being ``alien'' to Euler.

To be sure, rhetorical and formal aspects of historical mathematics
can be fruitfully studied in their own right.  Yet an overemphasis on
the rhetorical aspect to a point of dismissing as ``extraneous''
scholarly work that chooses to focus on the Eulerian mathematics
\emph{per se}, is untenable.

One may well wonder whether it sheds more light on Euler to observe,
as Laugwitz does, that Euler's infinitesimal procedures (Reeder's
\emph{inferential moves}) turn out to depend on \emph{hidden lemmas}
(such as those concerning estimates for infinite sums of
infinitesimals) but are otherwise remarkably robust and formalizable
in modern infinitesimal mathematics; or whether it sheds more light to
assert nonchalantly, as Ferraro does, that Euler considered
infinitesimals to be exactly zero as a kind of \emph{rhetorical
device}, and that therefore there are neither ``gaps" nor ``hidden
lemmas" in his proofs.  Relating to Euler's substitution
\[
\frac{i-1}{i}=1
\]
as a \emph{rhetorical device} as Ferraro does fails to explain why
Euler sometimes disallows this type of substitution, as when Euler
explains that
\[
1+\frac{x}{i}
\]
\emph{cannot} be replaced by~$1$ in factors of an infinite product in
the passage from \cite[\S\,156]{ai} cited in Section~\ref{31a}.
This passage from Euler explicitly contradicts Ferraro's
\emph{rhetorical} reading.

\subsection{Euler \emph{vs} Berkeley, H.~M.~Edwards, and Gray}
\label{101}

Cleric Berkeley's critique tends to receive exaggerated attention in
the literature.  We second Fraser's assessment to the effect that
``Berkeley's critique seems to have limited intrinsic merit''
\cite[p.\;453, note~3]{Fr}.  We now examine Ferraro's approach to this
critique.  Ferraro states that
\begin{quote}
a [sic] unproblematic translation of certain chapters in the history
of mathematics into modern terms tacitly assumes that the same logical
and conceptual framework guiding work in modern mathematics also
guided work in past mathematics \cite[p.\;2]{Fe12}.
\end{quote}
Here Ferraro expresses a legitimate concern.  Certainly one shouldn't
project the conceptual framework guiding modern mathematics, upon an
18th century text.  However, in the very next paragraph, Ferraro
proceeds to state: ``[Berkeley] did not cast any doubt upon the
usefulness of the calculus in solving many problems of physics or
geometry; nevertheless, he believed that it did not possess solid
\emph{foundations}'' (ibid.)  (emphasis added).  Let us now examine
the said foundations.

Berkeley's ``foundations,'' if any are to be found, amount to an
empiricist postulation of a minimal perceptual magnitude below which
one cannot descend, and a consequent rejection of an infinitely
divisible ``extension'' (i.e., continuum).  This is clearly not the
sense of the term \emph{mathematical foundations} that Ferraro has in
mind.  Rather than being concerned with the latter, Berkeley voiced
two separate criticisms: a \emph{metaphysical} and a \emph{logical}
one; see \cite{She87}.  The logical criticism concerns the alleged
inconsistency expressed by the conjunction~$(dx\not=0)\wedge(dx=0)$;
see Section~\ref{21}.  The metaphysical criticism is fueled by
Berkeley's empiricist doubts about entities that are below any finite
perceptual threshold.

Ferraro's description of Berkeley's criticism in terms of
``foundations'' falls prey to the very shortcoming he seeks to
criticize, namely grafting modern concepts upon ones exploited in
historical mathematics.  

A related attempt by H. M. Edwards to sweep Euler's infinitesimals
under an Archimedean rug in \cite{Ed07b} was analyzed in \cite{Ka15b}.
Edwards recently attempted to defend his comment that Euler's
infinitesimal computations
\begin{quote}
will not find a receptive audience today, when students are taught to
shrink from differentials as from an \emph{infectious disease}
\cite[p.~579]{Ed07b}, \cite[p.~52]{Ed16} (emphasis added)
\end{quote}
against our criticism in \cite{Ka15b}.  In recent years it has become
popular to interpret differentials as 1-forms.  This is fine, but it
is not Euler's view, as we show in \cite{Ka15b} and in the present
work.  In his response to that article, Edwards clarifies that he does
not dismiss Euler's use of differentials the way many others do.  But
it is not Edwards' \emph{disposition} toward differentials that is the
problem, but rather his \emph{interpretation} of Euler's
differentials.  In his response, Edwards again fails to acknowledge
that Euler's use of \emph{bona fide} infinitesimals is not reducible
to a purely algebraic algorithm.

Instead, Edwards indulges in rhetorical non-sequiturs against
Robinson's framework, accusing it of being ``far stranger than
anything Euler could have imagined.''  Edwards further accuses the
authors of \cite{Ka15b} of ``entertain[ing] \emph{strange ideas} about
the concept of the infinite'' (emphasis added).  However, Edwards'
remarks amount to a baseless ad hominem attack, since the article in
question said not a word about either Robinson or his framework,
focusing instead on the shortcomings of Edwards' take on Euler's work,
including a forced constructivist paraphrase thereof and an
anachronistic misattribution of the notion of \emph{derivative} to
Euler.

The book \cite{Ed79} (unrelated) presents a sympathetic view of
Robinson's framework, as does the book \cite{Ta14} which presents
ultraproducts as a bridge between discrete and continous analysis.

A year after the publication of H. Edwards' misguided analysis of
Euler in \cite{Ed07b}, J. Gray claimed that ``Euler's attempts at
explaining the foundations of calculus in terms of differentials,
which are and are not zero, are \emph{dreadfully weak}.''
\cite[p.\;6]{Gr08b} (emphasis added) Prisoner of A-track methodology,
Gray does not fail to succumb to Weierstrass's ghost when he claims in
his \emph{Plato's ghost} that Cauchy ``defined what it is for a
function \ldots to be continuous \ldots using careful, if not
altogether unambiguous, \emph{limiting arguments}.''
\cite[p.\;62]{Gr08a} (emphasis added) 

\emph{Pace} Gray, it is inaccurate to claim that Cauchy defined
continuity using \emph{limiting arguments}.  The word \emph{limit}
does appear in Cauchy's infinitesimal definition of continuity
(reproduced only two pages later in \emph{Plato's ghost}): ``the
function~$f(x)$ is continuous with respect to~$x$ between the given
limits if, between these limits, an infinitely small increment in the
variable always produces an infinitely small increment in the function
itself.'' \cite[p.\;26]{BS} Evidently, \emph{limits} do appear in
Cauchy's definition (though they are replaced by \emph{bounds} in
\cite[p.\;64]{Gr08a}).  However, they appear only in the sense of the
endpoints of the interval, rather than any sense related to the
Weierstrassian notion of the limit.

Gray's grafting of Weierstrassian limits upon Cauchy's definition of
continuity comes at a high price in anachronism.  For a recent study
of Cauchy based on Robinson's framework see \cite{CM}.

Ferraro could have used another term in place of \emph{foundations};
however, the exaggerated significance attached to Berkeley's allegedly
\emph{foundational} critique becomes apparent when Ferraro declares
that
\begin{quote}
The \emph{crux of the question} lay in knowing what meaning to
attribute to the equation~$a+dx=a$.  The exactness of mathematics
required, according to Euler, that the differential~$dx$ should be
precisely equal to~$0$: simply by assuming that~$dx=0$, the outrageous
attacks on the calculus would be shown to lack any basis
\cite[p.\;3]{Fe12} (emphasis added).
\end{quote}
Is this really the ``crux of the question'' as Ferraro contends?  As
discussed in Section~\ref{e36}, the \emph{exact zero} infinitesimals
are untenable and lead to insoluble paradoxes.  Meanwhile, the answer
to Berkeley's logical criticism lies elsewhere, namely the generalized
notion of equality implied by both the Leibnizian \emph{transcendental
law of homogeneity} and the Eulerian \emph{geometric comparison} (see
Section~\ref{22}) dubbed \emph{the principle of cancellation} by
Ferraro.  Characterizing Berkeley's logical criticism as the ``crux of
the question'' exaggerates the significance of his flawed empiricist
critique of infinitesimals.

\subsection{Aristotelian continuum?}
\label{e34}

Euler defined \emph{quantity} as that which could be increased or
reduced in his \emph{Elements of Algebra}: ``Whatever is capable of
increase or diminution, is called magnitude, or quantity \cite{Eu}.''

This may have been a common definition in Euler's time, but it was not
the classical definition of quantity.  What is called today the
Archimedean axiom characterizes the ancient Greek notion of quantity,
but it does not appear in modern mathematics until 1885 when it was
rediscovered in \cite{Sto85}.  Ferraro claims that
\begin{quote}
(1) Euler did not have the mathematical concept of set, nor the theory
of real numbers nor the modern notion of function.  (2) \emph{He based
  the calculus on the classic notion of quantity}.  (3) Quantity was
conceived of as that which could be increased or reduced \cite{Fe12}
(emphasis and numerals added).
\end{quote}
Ferraro's first and last claims are beyond dispute, but his
intermediate claim (italicized above) is dubious.  Namely, the claim
that Euler's notion of quantity was a ``classic" one is unsupported by
evidence.  Ferraro seeks to connect Euler's quantity to the notion of
quantity of unspecified \emph{ancient Greeks} as well as to \emph{the
classical Aristotelian conception}: ``[T]he Eulerian continuum is a
slightly modified version of the Leibnizian continuum, as described by
Breger [1992a, 76--84], which, in turn, has many aspects in common
with the classical Aristotelian conception.''  \cite[p.\;37]{Fe04}
Here Ferraro is referring to \cite{Br92}.  Breger does write on
page~76 that ``Leibniz reprend la th\'eorie aristot\'elicienne du
continu'' but in the same sentence he continues: ``en y apportant
trois modifications.''  One of these modifications, according to
Breger, is ``l'emploi des grandeurs infinit\'esimales.''  Ferraro's
claim that Breger's description of the Leibnizian continuum has ``many
aspects in common'' with the Aristotelian one appears to misrepresent
Breger's position as far as infinitesimals are concerned.

Thus, while the Archimedean axiom belongs to the \emph{classical} and
\emph{modern} notions of magnitude, it is found neither in
\emph{Euler's} characterisation of quantity as cited above, nor in
Leibniz's view of quantity.  See Sections~\ref{s38}, \ref{s39}, and
\ref{121b} for a discussion of quantity from Euclid to Euler.

\section{Conclusion}
\label{hacking}

In his essay for the collection \emph{Euler reconsidered}, Ferraro
writes: ``Euler's tripartite division of analysis was also the
manifestation of his aim to reduce analysis as far as possible to
algebraic notions; this latter term is used here to refer to notions
deriving from an \emph{infinitary extension of the principles of
analysis of finite quantities}'' \cite[p.\;45]{Fe07} (emphasis added).

\subsection{Cantor's ghost}

Ferraro's reference to Euler's \emph{infinitary extension of the
principles of analysis of finite quantities} alludes to concepts such
as infinite numbers and the associated infinite sums, or series, and
infinite products.  Infinite series and products are familiar
syntactic features of modern, A-track, analysis as formalized by
Cantor, Dedekind, and Weierstrass starting in the 1870s.  We would
like to comment on syntactic features that are noticeably absent from
the said analysis.  Cantor's own position can be briefly summarized as
follows:
\begin{quote}
Infinity, yes.

\noindent
Infinitesimals, no.
\end{quote}
In more detail, J. Dauben wrote:
\begin{quote}
Cantor devoted some of his most vituperative correspondence, as well
as a portion of the \emph{Beitr\"age}, to attacking what he described
at one point as the `infinitesimal Cholera bacillus of mathematics',
which had spread from Germany through the work of Thomae, du Bois
Reymond and Stolz, to infect Italian mathematics 
\cite[pp.\;216-217]{Da80}.
\end{quote}
Dauben continues:
\begin{quote}
Any acceptance of infinitesimals necessarily meant that his own theory
of number was incomplete.  Thus to accept the work of Thomae, du
Bois-Reymond, Stolz and Veronese was to deny the perfection of
Cantor's own creation.  Understandably, Cantor launched a thorough
campaign to discredit Veronese's work in every way possible (ibid.)
\end{quote}
Ferraro elaborates on his \emph{infinitary} comment cited above as
follows: ``\emph{Euler was not entirely successful in achieving his
aim, since he introduced infinitesimal considerations in various
proofs}; however, algebraic analysis, as a particular field of
mathematics, was clearly set out in the \emph{Introductio}.''
\cite[p.\;45]{Fe07} (emphasis added).  Given Ferraro's acknowledgment
that Euler exploits an \emph{infinitary extension of the principles of
analysis of finite quantities} as cited above, one might have expected
that such an infinitary extension involves both infinite numbers
\emph{and} infinitesimals.  

Yes Ferraro appears to feel, apparently following Cantor, that
infinite series constitute legitimate and \emph{successful} infinitary
extensions, whereas inferences involving infinitesimals do \emph{not}.
However, infinite numbers~$i$ and infinitesimals~$\omega$ in Euler are
related by the simple equation
\[
\omega=\frac{1}{i}
\]
or more generally~$i\omega=k$ where~$k$ is finite.  Why would one be
successful and the other \emph{not entirely successful}?  A possible
source of the distinction is the reliance on a conceptual framework
where infinite series admit suitable A-track proxies, whereas
infinitesimals do not.  Ferraro continues: ``At the end of the
eighteenth century, Euler's plan to undertake an algebraic treatment
of the broadest possible part of analysis of infinity had far-reaching
consequences when Lagrange tried to reduce the whole of calculus to
algebraic notions\ldots{}'' \cite[p.\;45]{Fe07}.  Such a \emph{broadest
possible} algebraic framework is apparently not broad enough, in
Ferraro's view, to encompass Eulerian infinitesimals.

\subsection{Primary point of reference?}

In a similar vein, Fraser claims that
\begin{quote}
\ldots{} classical analysis developed out of the older subject and it
remains a \emph{primary point of reference} for understanding the
eighteenth-century theories.  By contrast, nonstandard analysis and
other non-Archimedean versions of calculus emerged only fairly
recently in somewhat \emph{abstruse} mathematical settings that bear
\emph{little connection} to the historical developments one and a
half, two or three centuries earlier.  \cite[p. 27]{Fr15} (emphasis
added)
\end{quote}
For all his attempts to distance himself from Boyer's idolisation of
the triumvirate,%
\footnote{Historian Carl Boyer described Cantor, Dedekind, and
Weierstrass as \emph{the great triumvirate} in \cite[p.~298]{Boy}; the
term serves as a humorous characterisation of both A-track scholars
and their objects of adulation.}  
Fraser here commits himself to a position similar to Boyer's.  Namely,
Fraser claims that modern punctiform A-track analysis is \emph{a
primary point of reference} for understanding the analysis of the
past.  His sentiment that modern punctiform B-track analysis
\emph{bears little connection to the historical developments} reveals
insufficient attention to the procedure/ontology dichotomy.  A
sentiment of the inevitability of classical analysis is explicitly
expressed by Fraser who feels that ``classical analysis developed out
of the older subject and it remains a primary point of reference for
understanding the eighteenth-century theories'' yet his very
formulation involves circular reasoning.  It is only if one takes
classical analysis as \emph{a primary point of reference} that it
becomes plausible to conclude that it inevitably \emph{developed out
of the older subject}.

Such a position amounts to an unconditional adoption of the
teleological \emph{butterfly model} for the evolution of analysis,
where infinitesimals are seen as an evolutionary dead-end.
Elaborating an application of his butterfly/Latin dichotomy (see
Section~\ref{s32}) to the case of infinitesimals, Ian Hacking writes:
\begin{quote}
If analysis had stuck to infinitesimals in the face of philosophical
nay-sayers like Bishop Berkeley, analysis might have looked very
different.  Problems that were pressing late in the nineteenth
century, and which moved Cantor and his colleagues, might have
received a different emphasis, if any at all.  This alternative
mathematics might have seemed just as `successful', just as `rich', to
its inventors as ours does to us.  In that light, as Mancosu argued,
transfinite set theory now looks much more like the result of one of
Zeilberger's random walks than an inevitable mathematical
development. \cite[p.\;119]{Ha14}
\end{quote}

\subsection{Paradigm shift}

Laugwitz's pioneering articles from the 1980s such as \cite{Lau87},
\cite{La87}, and \cite{Lau89} built upon earlier studies, particularly
\cite{Ro66} and \cite{La}.  This work ushered in a new era in Euler
and Cauchy scholarship.  It became possible to dispense with, and go
beyond, the worn \emph{clich\'es} about unrigorous infinitesimalists
and their inconsistent manipulations with mystical infinitesimals.  In
the case of Euler, it became possible to formalize and interpret some
of his finest achievements in a way that sheds new light on the
\emph{methods} he used.  This work points to a coherence of his
formerly disparaged procedures based on the principle of cancellation,
infinitesimals, and infinite numbers, and establishes a historical
continuity in the procedures of infinitesimalists from Leibniz and
Euler to Robinson.

Such a paradigm shift in Euler scholarship has encountered resistance
from Ferraro, Fraser, Gray and other historians, who often cling to
Procrustean (and often slavishly post-Weierstrassian) frameworks of
Euler interpretation.  Thus, Gray finds Euler's explanations
``dreadfully weak'' but such a dismissive attitude toward Euler comes
at a high price in anachronism when applied to the 18th century.
Failing to distinguish clearly between procedural and ontological
issues, these historians focus on the latter and stress the obvious
point that modern set theory is alien to Euler's ontology, thus
falling back on strawman misrepresentions of the new wave of
scholarship.  The new scholarship accepts the obvious ontological
point, and focuses rather on the methodological issues of the
compatibility of Euler's \emph{inferential moves} and their proxies
provided by procedures available in modern infinitesimal frameworks.

Seeing with what dexterity Leibniz and Euler operated on infinite sums
as if they were finite sums, a modern scholar is faced with a stark
choice.  He can either declare that they didn't know the difference
between finite and infinite sums, or detect in their procedures a
unifying principle (explicit in the case of Leibniz, and more implicit
in the case of Euler) that, under suitable circumstances, allows one
to operate on infinite sums as on finite sums.  The former option is
followed by Ferraro, and is arguably dictated by self-imposed
limitations of an A-track interpretive framework.  The latter option
is the pioneering route of Robinson, Lakatos, Laugwitz, and others in
interpreting the infinitesimal mathematics of Leibniz, Euler, and
Cauchy.

\appendix

\section{Analysis of Euler's proof}
\label{31}

In Section~\ref{34} we summarized Euler's derivation of the product
decomposition for sine.  The derivation of infinite product
decompositions~\eqref{eu1} and \eqref{eu2} as found in
\cite[\S\,156]{ai} can be broken up into seven steps as follows.
Recall that Euler's~$i$ is an infinite integer.
\vspace{1ex}

\emph{Step 1}.  Euler observes that
\begin{eqnarray}
\label{eu3}
2\sinh x &=& e^x-e^{-x} =
\left(1+\frac x{i}\right)^{i}-\left(1-\frac x{i}\right)^{i},
\end{eqnarray}
where~${i}$ is an infinitely large natural number.  To motivate the
next step, note that the expression
$x^i-1=(x-1)(1+x+x^2+\ldots+x^{i-1})$ can be factored further as a
product $\prod_{k=0}^{i-1}(x-\zeta^k)$, where~$\zeta=e^{2\pi
\sqrt{-1}/i}$; conjugate factors can then be combined to yield a
decomposition into real quadratic terms.
\vspace{1ex}

\emph{Step 2}.  Euler uses the fact that~$a^{i}-b^{i}$ is the product
of the factors
\begin{equation}
\label{eu4}
a^2+b^2-2ab\cos\frac{2k\pi}{i}\,,\quad\text{where}\quad 1\le
k<\frac{i}{2}\,,
\end{equation}
together with the factor~$a-b$ and, if~${i}$ is an even number, the
factor~$a+b$, as well.
\vspace{1ex}

\emph{Step 3}.  Setting~$a=1+\frac x{i}$ and~$b=1-\frac x{i}$
in~\eqref{eu3}, Euler transforms expression~\eqref{eu4} into the form
\begin{equation}
\label{36b}
2+2\frac{x^2}{{i}^2}-2\bigg(1-\frac{x^2}{{i}^2}\bigg)
\cos\frac{2k\pi}{i}\,.
\end{equation}

\emph{Step 4}.  Euler then replaces~\eqref{36b} by the expression
\begin{equation}
\label{37}
\frac{4k^2\pi^2}{{i}^2}
\bigg(1+\frac{x^2}{k^2\pi^2}-\frac{x^2}{{i}^2}\bigg)\,,
\end{equation}
justifying this step by means of the formula
\begin{equation}
\label{38b}
\cos\frac{2k\pi}{i} = 1-\frac{2k^2\pi^2}{{i}^2}.
\end{equation}

\emph{Step 5}.  Next, Euler argues that the difference~$e^x-e^{-x}$ is
divisible by the expression
\begin{equation}
\label{39b}
1+\frac{x^2}{k^2\pi^2}-\frac{x^2}{i^2}
\end{equation}
from \eqref{37}, where ``we omit the term~$\frac{x^2}{{i}^2}$ since
even when multiplied by~${i}$, it remains infinitely small''
\cite{aie}.
\vspace{1ex}

\emph{Step 6}.  As there is still a factor of~$a-b=2x/{i}$, Euler
obtains the final equality~\eqref{eu1}, arguing that then ``the
resulting first term will be~$x$'' (in order to conform to the
Maclaurin series for~$\sinh x$).

\vspace{1ex} \emph{Step 7}.  Finally, formula~\eqref{eu2} is obtained
from \eqref{eu1} by means of the substitution~$x\mapsto
\sqrt{-1}\,x$.\qed\vspace{1ex}

Euler's argument in favor of \eqref{eu1} and \eqref{eu2} was
formalized in terms of a proof in Robinson's framework in \cite{lux}.
However, Luxemburg's formalisation deviates from Euler's argument
beginning with steps 3 and~4, and thus circumvents the most
problematic steps 5 and~6.  A proof in Robinson's framework,
formalizing Euler's argument step-by-step throughout, appeared in the
article \cite{kan}; see also \cite{mt} as well as the monograph
\cite[Section 2.4a]{kr}.  This formalisation interprets problematic
details of Euler's argument on the basis of general principles in
Robinson's framework, as well as general analytic facts that were
known in Euler's time.  Such principles and facts behind some early
proofs exploiting infinitesimals are sometimes referred to as
\emph{hidden lemmas} in this context; see \cite{Lau87}, \cite{Lau89},
and \cite{mt}.

For instance, a hidden lemma behind Step 4 asserts, on the basis of
the evaluation of the remainder~$R$ of the Taylor expansion
\[
\cos\frac{2k\pi}i=1-\frac{2k^2\pi^2}{i^2}+R\;,
\]
that the quadratic polynomial
$T_k(x)=2+2\frac{x^2}{{i}^2}-2\big(1-\frac{x^2}{{i}^2}\big)
\cos\frac{2k\pi}{i}$ as in~\eqref{36b} admits the representation
\[
T_k(x)= C_k\,\big(U_k(x)+p_k\cdot x^2\big)\,,
\]
where~$C_k$ and~$p_k$ do not depend on~$x$ while
\[
U_k(x)= 1+\frac{x^2}{k^2\pi^2}-\frac{x^2}{i^2}\,,
\]
and for any \emph{standard} real~$x$ and any finite or infinitely
large integer~$k\le\frac i2\,$ the following holds:
\begin{enumerate}
\item
if~$k$ is finite then~$p_k$ is infinitesimal, and\vspace{1ex}
\item
there is a real~$\gamma$ such that~$|p_k|<\gamma\cdot k^{-2}$ for any
infinitely large~$k\le\frac i2\,$.
\end{enumerate}
This allows one to infer that the effect of the transformation of
step~4 on the product of factors~\eqref{36b} is infinitesimal.  See
\cite[\S\,4]{kan} as well as equation~(11) on page~75 in \cite{kr} for
additional details.

Some \emph{hidden lemmas} of a different kind, related to basic
principles of nonstandard analysis, are discussed in
\cite[pp.\;43ff.]{mt} (see below).

What clearly stands out from Euler's argument is his explicit use of
infinitesimal quantities such as~\eqref{36b} and~\eqref{37}, as well
as the approximate formula~\eqref{38b} which holds ``up to'' an
infinitesimal of higher order.  Thus, Euler exploited \emph{bona fide}
infinitesimals, rather than merely ratios thereof, in a routine
fashion in some of his best work.

We now provide further technical details on a hyperreal interpretation
of Euler's proof of the product formula for the sine function.  Our
goal here is to indicate how Euler's inferential moves find modern
proxies in a hyperreal framework.

We discuss the \emph{hidden lemmas} related to basic principles of
nonstandard analysis following \cite[pp.\;43ff.]{mt}, where it is
argued that the Euler sine factorisation and similar constructions are
best understood in the context of the following \emph{hidden
definition} in terms of modern nonstandard analysis.  The following
definition is borrowed from \cite[p.\;44]{mt}.

\begin{quote}
\emph{Definition.}  A sum~$a_1+ a_2 + a_3 + \ldots$ is
Euler-convergent (E-convergent) if and only if
\begin{enumerate}
\item[(i)]~$a_k$ is defined by an elementary function,%
\footnote{The precise meaning of the modern term \emph{elementary
function} is discussed in \cite[p.\;43, footnote~23]{mt}.}
\item[(ii)] for all infinite%
\footnote {Here the terms \emph{finite} and \emph{infinite} correspond
to \emph{limited} and \emph{infinitely large} in the terminology of
\cite{mt}.}
$J$, the sum~$a_1+ a_2 + \ldots+a_J$ is finite, and 
\item[(iii)] for all infinite pairs~$J<K$, the sum~$a_J+ a_{J+1} +
\ldots+a_K$ is infinitesimal.
\end{enumerate}
Similarly, a product~$(1+b_l)(1+b_2)(1+b_3)\ldots$ is Euler-convergent
if and only if (i)~$b_k$ is defined by an elementary function,
(ii)~for all infinite~$J$, the product~$(1+b_1)(1+b_2)\ldots(1 + b_J)$
is finite, and (iii) for all infinitely large~$J<K$, the product~$(1 +
b_J)(1 + b_{J+1})\ldots(1 + b_K)$ differs infinitesimally from 1.
\end{quote}

Next, McKinzie and Tuckey present a series of \emph{hidden lemmas}
implicit in Euler's argument.  The first such \emph{hidden lemma}
asserts that if the sums~$a_1+ a_2 + \ldots$ and~$b_1+ b_2 + \ldots$
are E-convergent and~$a_k\simeq b_k$ (meaning that~$a_k-b_k$ is
infinitesimal) for all \emph{finite}~$k$, then
\[
a_1+ a_2 + \ldots+a_N\simeq b_1+ b_2 + \ldots+b_N
\]
for all~$N$ finite and infinite.  To prove this lemma, it suffices to
note that if~$a_k\simeq b_k$ holds for all finite~$k$, then, by
Robinson's lemma (see e.\,g., Theorem~2.2.12, p.\;62 in \cite{kr}),
there is an infinite~$K$ such that~$a_1+\dots+a_k\simeq
b_1+\ldots+b_k$ holds for all~$k\le K$.

The second \emph{hidden lemma} asserts a similar property for
products.  The third \emph{hidden lemma} asserts that if, for all
finite~$x$, the sums
\[
f(x)=a_0+a_1x+a_2x^2+\ldots \quad\text{and}\quad
g(x)=b_0+b_1x+b_2x^2+\ldots
\]
are E-convergent and we have [$f(x)\simeq g(x)$].  This means
that~$a_0+a_1x+a_2x^2+\ldots+a_Jx^J\simeq
b_0+b_1x+b_2x^2+\ldots+b_Kx^K$ for all infinite~$J,K$.  Note that the
choice of~$J,K$ is immaterial by (ii) and (iii) of the definition of
E-convergence.  Then~$a_n\simeq b_n$ for all~$n$ finite and infinite.
A detailed analysis in \cite{mt} shows that these three lemmas,
together with an additional sublemma, suffice to formalize Euler's
derivations step-by-step in a hyperreal framework.




\begin{thebibliography}{A}


\bibitem[Arthur 2008] {Ar08} Arthur, R.  ``Leery Bedfellows: Newton
and Leibniz on the Status of Infinitesimals.'' In \emph{Infinitesimal
Differences: Controversies between Leibniz and his Contemporaries},
ed. Ursula Goldenbaum and Douglas Jesseph, 7--30, Berlin and New York,
De Gruyter.


\bibitem[Bair et al.~2013]{Ba13} Bair, J.; B\l{}aszczyk, P.; Ely, R.;
Henry, V.; Kanovei, V.; Katz, K.; Katz, M.; Kutateladze, S.; McGaffey,
T.; Schaps, D.; Sherry, D.; Shnider,~S.  ``Is mathematical history
written by the victors?''  \emph{Notices of the American Mathematical
Society} \textbf{60}, no.~7, 886-904.

See \url{http://www.ams.org/notices/201307/rnoti-p886.pdf}

and \url{http://arxiv.org/abs/1306.5973}




\bibitem[Bascelli et al. 2014]{Ba14} Bascelli, T.; Bottazzi, E.;
Herzberg, F.; Kanovei, V.; Katz, K.; Katz, M.; Nowik, T.; Sherry, D.;
Shnider, S.  ``Fermat, Leibniz, Euler, and the gang: The true history
of the concepts of limit and shadow.''  \emph{Notices of the American
Mathematical Society} \textbf{61}, no.~8, 848-864.


\bibitem[Bascelli et al. 2016]{Ba16} Bascelli, T.; B\l{}aszczyk, P.;
Kanovei, V.; Katz, K.; Katz, M.; Schaps, D.; Sherry, D.  ``Leibniz vs
Ishiguro: Closing a quarter-century of syncategoremania.''
\emph{HOPOS: The Journal of the International Society for the History
of Philosophy of Science} \textbf{6}, no.~1, 117-147.  

See \url{http://dx.doi.org/10.1086/685645}

and \url{http://arxiv.org/abs/1603.07209}



\bibitem[Beckmann 1967/1968] {Bec} Beckmann, F.  ``Neue Gesichtspunkte
zum 5.  Buch Euklids.''  \emph{Archive for History Exact Sciences}
\textbf{4}, 1--144.


\bibitem[Bell 2008] {Be08} Bell, J.  \emph{A primer of infinitesimal
analysis. Second edition}.  Cambridge University Press, Cambridge.


\bibitem[Benacerraf 1965] {Be65} Benacerraf, P.  ``What numbers could
not be.''  \emph{Philos. Rev.}  \textbf{74}, 47--73.

\bibitem[Benci \& Di Nasso 2003]{BD} Benci, V.; Di Nasso, M.
``Numerosities of labelled sets: a new way of counting.''
\emph{Advances in Mathematics} \textbf{173}, no.~1, 50--67.


\bibitem[Berkeley 1734] {Be} Berkeley, G.  \emph{The Analyst, a
Discourse Addressed to an Infidel Mathematician.}



\bibitem[B\l{}aszczyk 2013] {Bl13} B\l{}aszczyk, P.  ``A note on Otto
H\"{o}lder's treatise \emph{Die Axiome der Quantit\"{a}t und die Lehre
vom Mass}.''  \emph{Annales Academiae Paedagogicae Cracoviensis.
Studia ad Didacticum Mathematicae} V, 129-142 [In Polish].
%
%


\bibitem[B\l aszczyk, Katz \& Sherry 2013]{BKS} B\l{}aszczyk, P.;
Katz, M.; Sherry, D.  ``Ten misconceptions from the history of
analysis and their debunking.''  \emph{Foundations of Science},
\textbf{18}, no.~1, 43--74.  See \url{http://arxiv.org/abs/1202.4153}

and \url{http://dx.doi.org/10.1007/s10699-012-9285-8}

%
%



\bibitem[B\l{}aszczyk et al.~2016]{Bl16} B\l{}aszczyk, P.; Kanovei,
V.; Katz, M.; Sherry, D.  ``Controversies in the Foundations of
Analysis: Comments on Schubring's Conflicts.''  \emph{Foundations of
Science}, online first.

See \url{http://dx.doi.org/10.1007/s10699-015-9473-4}


\bibitem[B\l{}aszczyk \& Mr\'owka 2013] {BM} B\l{}aszczyk, P.;
Mr\'owka, K.  \emph{Euklides, Elementy, Ksi{e}gi V-VI. T\l{}umaczenie
i komentarz} [Euclid, Elements, Books V-VI.  Translation and
commentary].  Copernicus Center Press, Krak\'ow.


\bibitem[Borovik \& Katz 2012]{BK} Borovik, A.; Katz, M.  ``Who gave
you the Cauchy--Weier\-strass tale?  The dual history of rigorous
calculus.''  \emph{Foundations of Science} \textbf{17}, no.~3,
245-276.  See \url{http://dx.doi.org/10.1007/s10699-011-9235-x}



\bibitem[Bos 1974] {Bos} Bos, H.  ``Differentials, higher-order
differentials and the derivative in the Leibnizian calculus.''
\emph{Archive for History of Exact Sciences} \textbf{14}, 1--90.


\bibitem[Bos 2010] {Bos12} Bos, H.  Private communication, 2 november
2010.


\bibitem[Boyer 1949]{Boy} Boyer, C.  \emph{The concepts of the
calculus}.  Hafner Publishing Company.


\bibitem[Bradley \& Sandifer 2009]{BS} Bradley, R.; Sandifer, C.
\emph{Cauchy's Cours d'analyse.  An annotated translation}.  Sources
and Studies in the History of Mathematics and Physical
Sciences. Springer, New York.

\bibitem[Breger 1992] {Br92} Breger, H.  ``Le continu chez Leibniz.''
In \emph{Le labyrinthe du continu} (Cerisy-la-Salle, 1990), 76-84,
Springer, Paris.


\bibitem[Carroll et al. 2013]{Ca13} Carroll, M.; Dougherty, S.;
Perkins, D.  ``Indivisibles, Infinitesimals and a Tale of
Seventeenth-Century Mathematics.''  \emph{Mathematics Magazine}
\textbf{86}, no.~4, 239--254.



\bibitem[Cauchy 1823] {Ca23} Cauchy, A. L.  \emph{R\'esum\'e des Le\c
cons donn\'ees \`a l'Ecole Royale Polytechnique sur le Calcul
Infinit\'esimal}.  Paris, Imprim\'erie Royale, 1823.  In \emph{Oeuvres
compl\`etes}, Series 2, Vol. 4, pp.\;9--261. Paris: Gauthier-Villars,
1899.


\bibitem[Ciesielski \& Miller 2016]{CM} Ciesielski, K.; Miller, D.
``A continuous tale on continuous and separately continuous
functions.''  \emph{Real Analysis Exchange} \textbf{41}, no.~1,
19--54.

\bibitem[Clavius 1589]{Clavius} Clavius, C.  \emph{Euclidis
Elementorum. Libri XV}, Roma.

\bibitem[Dauben 1980] {Da80} Dauben, J.  ``The development of
Cantorian set theory.''  In \emph{From the calculus to set theory,
1630-1910}, 181-219, Princeton Paperbacks, Princeton University Press,
Princeton, NJ, 2000. [Originally published in 1980.]


\bibitem[De Risi 2016]{De16} De Risi, V.  ``The Development of
Euclidean Axiomatics.  The systems of principles and the foundations
of mathematics in editions of the Elements from Antiquity to the
Eighteenth Century.''  \emph{Archive for History of Exact Science},
online first.  See \url{http://dx.doi.org/10.1007/s00407-015-0173-9}



\bibitem[Di Nasso \& Forti 2010]{DF} Di Nasso, M.; Forti, M.
``Numerosities of point sets over the real line.''  \emph{Transactions
of the American Mathematical Society} \textbf{362}, no.~10,
5355--5371.


\bibitem[Dijksterhuis 1987] {Di87} Dijksterhuis, D.  \emph{Archimedes
1987}.  translated by C. Diksboorn, Princeton Univ Press, Princeton
(first printed in 1956 by E. Muskgaard), pp.\;130-133.
%


\bibitem[Edwards, C. H., Jr. 1979]{Ed79} Edwards, C. H., Jr.
\emph{The historical development of the calculus}.  Springer-Verlag,
New York--Heidelberg.

\bibitem[Edwards 2007] {Ed07b} Edwards, H.  ``Euler's definition of
the derivative.''  \emph{Bull. Amer. Math. Soc. (N.S.)} \textbf{44},
no.~4, 575--580.




\bibitem[Edwards 2015] {Ed16} Edwards, H.  ``Euler's conception of the
derivative.''  \emph{The Mathematical Intelligencer} \textbf{47},
no~4, 52--53.

See \url{http://dx.doi.org/10.1007/s00283-015-9560-y}




\bibitem[Ehrlich 2006] {Eh06} Ehrlich, P.  ``The rise of
non-Archimedean mathematics and the roots of a misconception.  I. The
emergence of non-Archimedean systems of magnitudes.''  \emph{Archive
for History of Exact Sciences} \textbf{60}, no.~1, 1--121.


\bibitem[Euclid 1660]{Eu60} Euclid.  \emph{Euclide's Elements; The
whole Fifteen Books, compendiously Demonstrated}.  By Mr.~Isaac Barrow
Fellow of \emph{Trinity College} in Cambridge.  And Translated out of
the Latin.  London.


\bibitem[Euclid 2007]{Euc} Euclid.  \emph{Euclid's Elements of
Geometry}.  Edited, and provided with a modern English translation, by
Richard Fitzpatrick.  See

\url{http://farside.ph.utexas.edu/euclid.html}


\bibitem[Euler 1730-1731] {Eu30} Euler, L.  ``De progressionibus
trascendentibus seu quarum termini generales algebraice dari
nequeunt.''  In: Euler, \emph{Opera omnia} (Series I, Opera
Mathematica, Berlin, Bern, Leipzig, 1911--\ldots) 14 (1), 1-24.



\bibitem[Euler 1748] {ai} Euler, L.  \emph{Introductio in Analysin
Infinitorum, Tomus primus}.  Saint Petersburg and Lausana.


\bibitem[Euler 1755] {ed} Euler, L.  \emph{Institutiones Calculi
Differentialis}.  Saint Petersburg.


\bibitem[Euler 1988] {aie} Euler, L.  \emph{Introduction to analysis
of the infinite.  Book I}.  Translated from the Latin and with an
introduction by J.~Blanton.  Springer-Verlag, New York.  (translation
of \cite{ai})



\bibitem[Euler 2000] {Eu55} Euler, L.  \emph{Foundations of
Differential Calculus}.  English translation of Chapters 1--9 of
\cite{ed} by J.\,Blanton, Springer, N.Y.



\bibitem[Euler 1768-1770] {Eu68} Euler, L.  \emph{Institutionum
calculi integralis}.



\bibitem[Euler 1771] {Eu71} Euler, L.  \emph{Vollst\"{a}ndige
Anleitung zur Algebra}.  Kaiserliche Akademie der Wissenschaften,
St.~Petersburg.

\bibitem[Euler 1807] {Eu07} Euler, L.  \emph{\'El\'ements
d'alg\`ebre}.  Courcier, Paris.

\bibitem[Euler 1810] {Eu} Euler L.  \emph{Elements of Algebra}.
Translated form the French with additions of La Grange, Johson and
Co., London.


\bibitem[Ferraro 1998] {Fe98} Ferraro, G.  ``Some aspects of Euler's
theory of series: inexplicable functions and the Euler-Maclaurin
summation formula.''  \emph{Historia Mathematica} \textbf{25}, no.~3,
290--317.


\bibitem[Ferraro 2004] {Fe04} Ferraro, G.  ``Differentials and
differential coefficients in the Eulerian foundations of the
calculus.''  \emph{Historia Mathematica} \textbf{31}, no.~1, 34--61.


\bibitem[Ferraro 2007] {Fe07} Ferraro, G.  ``Euler's treatises on
infinitesimal analysis: \emph{Introductio in analysin infinitorum,
institutiones calculi differentialis, institutionum calculi
integralis}.''  In \emph{Euler reconsidered}, 39-101, Kendrick Press,
Heber City, UT.


\bibitem[Ferraro 2008]{Fe08} Ferraro, G.  \emph{The rise and
development of the theory of series up to the early 1820s.}  Sources
and Studies in the History of Mathematics and Physical Sciences.
Springer, New York.

\bibitem[Ferraro 2012] {Fe12} Ferraro, G.  ``Euler, infinitesimals and
limits.''  Manuscript (february~2012).  See
\url{http://halshs.archives-ouvertes.fr/halshs-00657694}


\bibitem[Ferraro \& Panza 2003] {FP} Ferraro, G., Panza, M.
``Developing into series and returning from series: a note on the
foundations of eighteenth-century analysis.''  \emph{Historia
Mathematica} \textbf{30}, no.~1, 17--46.


\bibitem[Fraenkel 1928] {Fran} Fraenkel, A.  \emph{Einleitung in die
Mengenlehre}.  Dover Publications, New York, N. Y., 1946 [originally
published by Springer, Berlin, 1928].

\bibitem[Fraser 1999] {Fr} Fraser, C.  ``Book review: Paolo Mancosu.
Philosophy of mathematics and mathematical practice in the seventeenth
century.  The Clarendon Press, Oxford University Press, New York.''
\emph{Notre Dame Journal of Formal Logic} \textbf{40}, no.~3,
447--454.

\bibitem[Fraser 2015]{Fr15} Fraser, C.  ``Nonstandard analysis,
infinitesimals, and the history of calculus.''  In \emph{A Delicate
Balance: Global Perspectives on Innovation and Tradition in the
History of Mathematics}, D. Row \& W. Horng, eds., Birkh\"auser,
Springer, pp.~25--49.

\bibitem[Gerhardt 1850-1863] {Ge50} Gerhardt, C. I. (ed.)
\emph{Leibnizens mathematische Schriften}.  Berlin and Halle: Eidmann.
%
%


\bibitem[G\"odel 1990] {Go90} G\"odel, K.  \emph{Collected Works}.
Feferman, S., Dawson, J., Kleene, S., Moore, G., Solovay, R., and van
Heijenoort, J., editors.  Vol. II.  New York: Oxford University Press.


\bibitem[Goldblatt 1998] {Go} Goldblatt, R.  \emph{Lectures on the
hyperreals.  An introduction to nonstandard analysis}.  Graduate Texts
in Mathematics \textbf{188}.  Springer-Verlag, New York.


\bibitem[Gordon, Kusraev \& Kutateladze 2002] {GKK} Gordon, E.,
Kusraev, A., Kutateladze, S.  \emph{Infinitesimal analysis}.  Updated
and revised translation of the 2001 Russian original.  Translated by
Kutateladze.  \emph{Mathematics and its Applications} \textbf{544}.
Kluwer Academic Publishers, Dordrecht.


\bibitem[Grant 1974]{Grant} Grant, E.  ``The definitions of Book V of
Euclid's Elements in thirteenth-century version, and
commentary. Campanus of Novara.''  In E.~Grant, \emph{A Source book in
medieval science}, Harvard University Press, 136--149.

\bibitem[Gray 2008a]{Gr08a} Gray, J.  \emph{Plato's ghost.  The
modernist transformation of mathematics}.  Princeton University Press,
Princeton, NJ.

\bibitem[Gray 2008b]{Gr08b} Gray, J.  ``A short life of Euler.''
\emph{BSHM Bull}.  \textbf{23}, no.~1, 1--12.


\bibitem[Hacking 2014] {Ha14} Hacking, I.  \emph{Why is there
philosophy of mathematics at all?}  Cambridge University Press.

\bibitem[Heiberg 1881] {Heib81} Heiberg, J.  \emph{Archimedis Opera
Omnia cum Archimedis Opera Omnia cum Commentariis Eutocii, vol. I}.
Teubner, Leipzig.

\bibitem[Heiberg 1883-1888] {Heib83} Heiberg, J.  \emph{Euclidis
Elementa, vol. I-V}.  Teubner, Leipzig.

\bibitem[H\"{o}lder 1901] {Ho} H\"{o}lder, O.  ``Die Axiome der
Quantit\"{a}t und die Lehre vom Mass.''  Berichte \"{u}ber die
Verhandlungen der K\"{o}niglich S\"{a}chsischen Gesellschaft der
Wissenschaften zu Leipzig.  Mathematisch-Physische Classe, 53,
Leipzig, 1--63.
%
%

\bibitem[H\"{o}lder 1996] {Ho96} H\"{o}lder, O.  ``The axioms of
quantity and the theory of measurement.''  Translated from the 1901
German original and with notes by Joel Michell and Catherine
Ernst. With an introduction by Michell.  \emph{J. Math. Psych.}
\textbf{40}, no. 3, 235--252.



\bibitem[Ishiguro 1990] {Is} Ishiguro, H.  \emph{Leibniz's philosophy
of logic and language.  Second edition}.  Cambridge University Press,
Cambridge.


\bibitem[Jesseph 2015]{Je15} Jesseph, D.  ``Leibniz on the Elimination
of infinitesimals.''  In \emph{G.W.Leibniz, Interrelations between
Mathematics and Philosophy}, pp.\,189--205.  Norma B. Goethe, Philip
Beeley, and David Rabouin, eds.  Archimedes Series \textbf{41},
Springer Verlag.


\bibitem[Kanovei 1988] {kan} Kanovei, V.  ``The correctness of Euler's
method for the factorization of the sine function into an infinite
product.''  \emph{Russian Mathematical Surveys} \textbf{43}, 65--94.


\bibitem[Kanovei 1993]{Ka93} Kanovei, V.  ``Leonhard Euler's summation
of a series of factorials that are alternating in sign.''  (Russian)
\emph{Istor.-Mat.~Issled.}  \textbf{34}, 16--46.


\bibitem[Kanovei, Katz \& Nowik 2016] {KKN} Kanovei, V.; Katz, K.;
Katz, M.; Nowik, T.  ``Small oscillations of the pendulum, Euler's
method, and adequality.''  \emph{Quantum Studies: Mathematics and
Foundations}.

See \url{http://dx.doi.org/10.1007/s40509-016-0074-x}

and \url{http://arxiv.org/abs/1604.06663}



\bibitem[Kanovei, Katz \& Mormann 2013] {KKM} Kanovei, V.; Katz, M.;
Mormann, T.  ``Tools, Objects, and Chimeras: Connes on the Role of
Hyperreals in Mathematics.''  \emph{Foundations of Science}
\textbf{18}, no.~2, 259--296.  

See \url{http://dx.doi.org/10.1007/s10699-012-9316-5}

and \url{http://arxiv.org/abs/1211.0244}


\bibitem[Kanovei, Katz \& Schaps 2015] {Ka15a} Kanovei, V., Katz, K.,
Katz, M., Schaps,~M.  ``Proofs and retributions, Or: Why Sarah can't
\emph{take} limits.''  \emph{Foundations of Science} \textbf{20},
no.~1, 1--25.  

See \url{http://dx.doi.org/10.1007/s10699-013-9340-0}


\bibitem[Kanovei, Katz \& Sherry 2015]{Ka15b} Kanovei, V.; Katz, K.;
Katz, M.; Sherry, D.  ``Euler's lute and Edwards' oud."  \emph{The
Mathematical Intelligencer} \textbf{37}, no.~4, 48--51.

See \url{http://dx.doi.org/10.1007/s00283-015-9565-6}

and \url{http://arxiv.org/abs/1506.02586}




\bibitem[Kanovei \& Reeken 2004] {kr} Kanovei, V., Reeken, M.
\emph{Nonstandard analysis, axiomatically}.  Springer Monographs in
Mathematics, Berlin, Springer.

\bibitem[Katz \& Katz 2011]{KK11} Katz, K.; Katz, M. ``Cauchy's
continuum."  \emph{Perspectives on Science} \textbf{19}, no.~4,
426-452.

See \url{http://dx.doi.org/10.1162/POSC_a_00047}

and \url{http://arxiv.org/abs/1108.4201}


\bibitem[Katz, Schaps \& Shnider 2013]{KSS13} Katz, M.; Schaps, D.;
Shnider, S.  Almost Equal: The Method of Adequality from Diophantus to
Fermat and Beyond.  \emph{Perspectives on Science} \textbf{21} (2013),
no.~3, 283-324.

See \url{http://www.mitpressjournals.org/doi/abs/10.1162/POSC_a_00101}

and \url{http://arxiv.org/abs/1210.7750}



\bibitem[Katz \& Sherry 2012] {KS2} Katz, M., Sherry, D.  ``Leibniz's
laws of continuity and homogeneity.''  \emph{Notices of the American
Mathematical Society} \textbf{59}, no.~11, 1550-1558.

See \url{http://www.ams.org/notices/201211/rtx121101550p.pdf}

and \url{http://arxiv.org/abs/1211.7188}


\bibitem[Katz \& Sherry 2013] {KS1} Katz, M., Sherry, D.  ``Leibniz's
infinitesimals: Their fictionality, their modern implementations, and
their foes from Berkeley to Russell and beyond.''  \emph{Erkenntnis}
\textbf{78}, no.~3, 571--625.  

See \url{http://dx.doi.org/10.1007/s10670-012-9370-y}

and \url{http://www.ams.org/mathscinet-getitem?mr=3053644}

and \url{http://arxiv.org/abs/1205.0174}




\bibitem[Katz 2014] {Ka14} Katz, V.  Review of ``Bair et al., Is
mathematical history written by the victors?  \emph{Notices
Amer. Math. Soc.}  \textbf{60} (2013), no.~7, 886--904.''  

See \url{http://www.ams.org/mathscinet-getitem?mr=3086638}


\bibitem[Keisler 1986]{Ke86} Keisler, H. J.  \emph{Elementary
Calculus: An Infinitesimal Approach.}  Second Edition.  Prindle, Weber
\& Schimidt, Boston.  See online version at
\url{http://www.math.wisc.edu/~keisler/calc.html}


\bibitem[Klein 1908] {Kl08} Klein, F.  \emph{Elementary Mathematics
from an Advanced Standpoint.  Vol.~I.  Arithmetic, Algebra, Analysis.}
Translation by E. R. Hedrick and C. A. Noble [Macmillan, New York,
1932] from the third German edition [Springer, Berlin, 1924].
Originally published as \emph{Elementarmathematik vom h\"oheren
Standpunkte aus} (Leipzig, 1908).


\bibitem[Kock 2006] {Ko06} Kock, A.  \emph{Synthetic differential
geometry.  Second edition}.  London Mathematical Society Lecture Note
Series, 333.  Cambridge University Press, Cambridge.

\bibitem[Kuhn 1962] {Ku} Kuhn, T.  \emph{The structure of scientific
revolutions}.  University of Chicago Press.


\bibitem[Lakatos 1978] {La} Lakatos, I.  ``Cauchy and the continuum:
the significance of nonstandard analysis for the history and
philosophy of mathematics.''  \emph{Mathematical Intelligencer}
\textbf{1}, no. 3, 151--161 (essay originally appeared in 1966).


\bibitem[Laugwitz 1987a] {Lau87} Laugwitz, D.  ``Hidden lemmas in the
early history of infinite series.''  \emph{Aequationes Mathematicae}
\textbf{34}, 264--276.

\bibitem[Laugwitz 1987b]{La87} Laugwitz, D.  ``Infinitely small
quantities in Cauchy's textbooks.''  \emph{Historia Mathematica}
\textbf{14}, 258--274.


\bibitem[Laugwitz 1989] {Lau89} Laugwitz, D.  ``Definite values of
infinite sums: aspects of the foundations of infinitesimal analysis
around 1820.''  \emph{Archive for History of Exact Sciences}
\textbf{39}, no.~3, 195--245.


\bibitem[Laugwitz 1992a] {Lau92} Laugwitz, D.  ``Early delta functions
and the use of infinitesimals in research.''  \emph{Revue d'histoire
des sciences} \textbf{45}, no.~1, 115--128.

\bibitem[Laugwitz 1992b] {Lau92b} Laugwitz, D.  ``Leibniz' principle
and omega calculus.''  In \emph{Le labyrinthe du continu}
(Cerisy-la-Salle, 1990), 144--154, Springer, Paris.


\bibitem[Laugwitz 1999]{La99} Laugwitz, D.  \emph{Bernhard Riemann
1826-1866.  Turning points in the conception of mathematics}.
Translated from the 1996 German original by Abe Shenitzer with the
editorial assistance of the author, Hardy Grant, and Sarah Shenitzer.
Birkh\"auser Boston, Boston.


\bibitem[Leibniz 1684] {Le84} Leibniz, G.  ``Nova methodus pro maximis
et minimis \ldots'' in \emph{Acta Erud.}, Oct. 1684.  See \cite{Ge50},
V, pp.~220--226.



\bibitem[Leibniz 1695] {Le95a} Leibniz, G.  To l'Hospital, 21 june
1695, in \cite{Ge50}, I, pp.~287--289.


\bibitem[Leibniz 1702]{Le02} Leibniz, G.  To Varignon, 2 febr., 1702,
in \cite{Ge50} IV, pp.~91--95.


\bibitem[Leibniz 1710] {Le10b} Leibniz, G.  ``Symbolismus memorabilis
calculi algebraici et infinitesimalis in comparatione potentiarum et
differentiarum, et de lege homogeneorum transcendentali.''  In
\cite[vol.~V, pp.\;377-382]{Ge50}.


\bibitem[l'H\^opital 1696] {lh} l'H\^opital, G.  \emph{Analyse des
Infiniment Petits pour l'Intelligence des Lignes Courbes}.


\bibitem[Luxemburg 1973] {lux} Luxemburg, W.  ``What is nonstandard
analysis?''  Papers in the foundations of mathematics.  \emph{American
Mathematical Monthly} \textbf{80}, no.~6, part II, 38--67.

\bibitem[Mancosu 1996] {Ma96} Mancosu, P.  \emph{Philosophy of
mathematics and mathematical practice in the seventeenth century.}
The Clarendon Press, Oxford University Press, New York.

\bibitem[Mancosu 2009] {Ma09} Mancosu, P.  ``Measuring the size of
infinite collections of natural numbers: was Cantor's theory of
infinite number inevitable?''  \emph{Rev. Symb. Log.} \textbf{2},
no.~4, 612--646.


\bibitem[McKinzie \& Tuckey 1997] {mt} McKinzie, M., Tuckey, C.
``Hidden lemmas in Euler's summation of the reciprocals of the
squares.''  \emph{Archive for History of Exact Sciences} \textbf{51},
29--57.



\bibitem[Mueller 1981] {Mue} Mueller, I.  \emph{Philosophy of
mathematics and deductive structure in Euclid's Elements}.  MIT Press,
Cambridge, Mass.--London [reprinted by Dover in 2006].


\bibitem[Nelson 1977]{Ne77} Nelson, E.  ``Internal set theory: a new
approach to nonstandard analysis.''  \emph{Bulletin of the American
Mathematical Society} \textbf{83}, no.~6, 1165--1198.


\bibitem[Nieuwentijt 1695] {Ni} Nieuwentijt, B.  \emph{Analysis
infinitorum, seu curvilineorum proprietates ex polygonorum natura
deductae.}  Amsterdam.

\bibitem[Nowik \& Katz 2015]{NK} Nowik, T., Katz, M.  ``Differential
geometry via infinitesimal displacements.''  \emph{Journal of Logic
and Analysis} \textbf{7}:5, 1--44.  See

\url{http://www.logicandanalysis.org/index.php/jla/article/view/237/106}
and \url{http://arxiv.org/abs/1405.0984}



\bibitem[Panza 2007] {Pa07} Panza, M.  ``Euler's \emph{Introductio in
Analysin Infinitorum} and the program of algebraic analysis.''  In
R. Backer (ed.)  \emph{Euler reconsidered}, 119--166, Kendrick Press.



\bibitem[Pulte 1998]{Pu98} Pulte, H.  ``Jacobi's criticism of
Lagrange: the changing role of mathematics in the foundations of
classical mechanics.''  \emph{Historia Mathematica} \textbf{25},
no.~2, 154--184.

\bibitem[Pulte 2012]{Pu12} Pulte, H.  ``Rational mechanics in the
eighteenth century.  On structural developments of a mathematical
science.''  \emph{Berichte zur Wissenschaftsgeschichte} \textbf{35},
no.~3, 183--199.


\bibitem[Quine 1968] {Qu} Quine, W.  ``Ontological Relativity.''
\emph{The Journal of Philosophy} \textbf{65}, no.~7, 185-212.


\bibitem[Reeder 2012] {Re12} Reeder, P.  ``A `Non-standard Analysis'
of Euler's \emph{Introductio in Analysin Infitorum}.''  MWPMW 13.
University of Notre Dame, october 27-28, 2012.  

See \url{http://philosophy.nd.edu/assets/81379/mwpmw_13.summaries.pdf}



\bibitem[Reeder 2013] {Re13} Reeder, P.  \emph{Internal Set Theory and
Euler's \emph{Introductio in Analysin Infinitorum}}.  MSc Thesis, Ohio
State University.


\bibitem[Robinson 1961] {Ro61} Robinson, A.  ``Non-standard
analysis.''  \emph{Nederl. Akad. Wetensch. Proc. Ser. A} \textbf{64} =
\emph{Indag. Math.} \textbf{23}, 432--440 (reprinted in Selected
Works, see \cite[pp.\;3--11]{Ro79}).


\bibitem[Robinson 1966] {Ro66} Robinson, A.  \emph{Non-standard
analysis}.  North-Holland Publishing Co., Amsterdam.


\bibitem[Robinson 1969] {Ro69} Robinson, A.  ``From a formalist's
points of view.''  \emph{Dialectica} \textbf{23}, 45--49.


\bibitem[Robinson 1979] {Ro79} Robinson, A.  \emph{Selected papers of
Abraham Robinson.  Vol. II.  Nonstandard analysis and philosophy}.
Edited and with introductions by W. A. J. Luxemburg and S. K\"orner.
Yale University Press, New Haven, Conn.


\bibitem[Sandifer 2007] {Sa} Sandifer, C.  ``Euler's solution of the
Basel problem--the longer story.''  In \emph{Euler at 300}, 105--117,
MAA Spectrum, Math. Assoc. America, Washington, DC.

\bibitem[Schubring 2016]{Sc15} Schubring, G.  ``Comments on a Paper on
Alleged Misconceptions Regarding the History of Analysis: Who Has
Misconceptions?'' \emph{Foundations of Science}, online first.  See
\url{http://dx.doi.org/10.1007/s10699-015-9424-0}

\bibitem[Sherry 1987] {She87} Sherry, D.  ``The wake of Berkeley's
Analyst: \emph{rigor mathematicae}?''  \emph{Stud. Hist. Philos. Sci.}
\textbf{18}, no.~4, 455--480.


\bibitem[Sherry \& Katz 2014] {SK} Sherry, D., Katz, M.
``Infinitesimals, imaginaries, ideals, and fictions.''  \emph{Studia
Leibnitiana} \textbf{44} (2012), no.~2, 166-192 (the article was
published in 2014 even though the journal issue lists the year as
2012). 

See \url{http://arxiv.org/abs/1304.2137}




\bibitem[Stolz 1883] {Sto83} Stolz, O.  ``Zur {G}eometrie der {A}lten,
insbesondere {\"u}ber ein {A}xiom des {A}rchimedes.''
\emph{Mathematische Annalen} \textbf{22} (4), 504--519.

\bibitem[Stolz 1885] {Sto85} Stolz, O.  \emph{Vorlesungen \"{u}ber
Allgemeine Arithmetik}, Teubner, Leipzig.


\bibitem[Tao 2014]{Ta14} Tao, T.  \emph{Hilbert's fifth problem and
related topics}.  Graduate Studies in Mathematics, 153. American
Mathematical Society, Providence, RI.

\bibitem[Tho 2012] {Th12} Tho, T.  ``Equivocation in the foundations
of Leibniz's infinitesimal fictions.''  \emph{Society and Politics}
\textbf{6}, no.~2, 70--98.

\bibitem[Vermij 1989] {Ver} Vermij, R.  ``Bernard Nieuwentijt and the
Leibnizian calculus.''  \emph{Studia Leibnitiana} \textbf{21}, no.~1,
69--86.


\bibitem[Wallis 2004] {Wa} Wallis, J.  (1656) \emph{The arithmetic of
infinitesimals}.  Translated from the Latin and with an introduction
by Jaequeline A. Stedall.  \emph{Sources and Studies in the History of
Mathematics and Physical Sciences}.  Springer-Verlag, New York.

\bibitem[Wallis 1685] {Wa85} Wallis, J.  \emph{Treatise of Algebra}.
Oxford.

\bibitem[Wartofsky 1976]{Wa76} Wartofsky, M.  ``The Relation Between
Philosophy of Science and History of Science.''  In R. S. Cohen,
P. K. Feyerabend, and M. W. Wartofsky (eds.), Essays in Memory of Imre
Lakatos, 717--737, Boston Studies in the Philosophy of Science XXXIX,
D. Reidel Publishing, Dordrecht, Holland.


\bibitem[Weber 1895] {We} Weber, H.  \emph{Lehrbuch der Algebra}.
Einleitung.  Vieweg, Braunschweig.



\end{thebibliography}
\end{document}